\documentclass[a4paper,10pt]{article}
\usepackage[pagewise]{lineno}
\usepackage{microtype}
\usepackage[english]{babel}
\usepackage[utf8]{inputenc}
\usepackage{hyperref}
\usepackage{bm}
\usepackage{sidecap}
\usepackage{float}
\usepackage{booktabs}
\usepackage{turnstile}
\usepackage{amssymb,amsthm,amsmath}
\usepackage{indentfirst}
\usepackage{color}
\usepackage{graphicx, caption, subcaption}
\usepackage{listingsutf8}
\usepackage{xcolor}
\usepackage{bbm}
\usepackage{url}
\usepackage{cite}
\usepackage{algorithm}
\usepackage{algorithmic}
\usepackage{bbm}
\usepackage{lipsum}
\usepackage{amsthm}
\usepackage{amsmath}

\newtheorem{remark}{Remark}
\newtheorem{definition}{Definition}
\newtheorem{theorem}{Theorem}[section]
\newtheorem{lemma}[theorem]{Lemma}

\newtheorem{property}{Property}

\usepackage[a4paper,top=3cm,bottom=3cm,left=3cm,right=3cm]{geometry}
\usepackage{verbatim}

\usepackage{chngcntr}
\newcommand{\ka}{\kappa}

\usepackage{cleveref}

\makeatletter
\def\BState{\State\hskip-\ALG@thistlm}
\makeatother

\usepackage{authblk}
\numberwithin{equation}{section}

\colorlet{corr}{magenta}

\usepackage{etoolbox} 

\title{Micro-Macro Decomposition of Particle Swarm Optimization Methods}

\author{Michael Herty}
\affil{\normalsize RWTH Aachen University, Institute of Geometry and Applied Mathematics, Templergraben 55, 52062 Aachen, Germany}
\affil{\normalsize Extraordinary Professor, University of Pretoria, Department of Mathematics and Applied Mathematics,  Private Bag X20, Hatfield 0028, South Africa\\
\texttt{herty@igpm.rwth-aachen.de}}

\author{Sara Veneruso\thanks{Corresponding author}}
\affil{\normalsize RWTH Aachen University, Institute of Geometry and Applied Mathematics, Templergraben 55, 52062, Aachen, Germany}
\affil{\normalsize University of Ferrara, Department of Mathematics and Computer Science, Via Machiavelli 30, 44121 Ferrara, Italy \\
\texttt{veneruso@igpm.rwth-aachen.de, sara.veneruso@unife.it}}

\date{August 18, 2025}

\begin{document}
\maketitle

\begin{abstract}
Solving non-convex minimization problems using multi-particle metaheuristic derivative-free optimization methods is still an active area of research. Popular methods are Particle Swarm Optimization (PSO) methods, that iteratively update a population of particles according to dynamics inspired by social interactions between individuals. We present a modification to include constrained minimization problems using exact penalization.
Additionally, we utilize the hierarchical structure of PSO to introduce a micro-macro decomposition of the algorithm. The probability density of particles is written as a convex combination of microscopic and macroscopic contributions, and both parts are propagated separately. The decomposition is dynamically updated based on heuristic considerations. Numerical examples compare the results obtained
using the algorithm in the microscopic scale, in the macroscopic scale, and using the new micro-macro decomposition.
\end{abstract}

\section{Introduction}
Many mathematical models nowadays study the coordinated movement of groups using multi-agent approaches \cite{vicsek,cucker}. Metaheuristic methods are based on simple rules of interaction between individual particles, which give complex results of self-organized and swarm intelligence \cite{albi,sumpter,motsch}. Some examples of these methods are Ant Colony Optimization \cite{dorigo, DORIGO2}, Genetic Algorithms \cite{holland, goldberg}, Simulated Annealing \cite{aarts,kirkpatrick}, Consensus-Based Optimization (CBO)\cite{Pinnau, carrillo,seung} and Particle Swarm Optimization (PSO) \cite{kennedy,kennedy2,shi, huang}. There are strong links between these methods, for example the CBO can be derived from the PSO via zero inertia limit \cite{cipriani}.

In this paper, we focus on the study of Particle Swarm Optimization, proposed for solving constrained minimization problems. The goal is to find 
\begin{equation}
\label{prob_init}
    \min_{x\in\mathbb{R}^d}\mathcal{F}(x) \quad \text{ subject to } x\in\mathcal{M}
\end{equation}
$d\ge 1$, where $\mathcal{F}(x)$ is a given continuous but possibly non convex objective function. It is assumed that $\mathcal{F}(x)$ admits an unique global minimizer $x^*_\mathcal{M}$ in $\mathcal{M}\subset \mathbb{R}^d$.

The PSO method solves the optimization problem \eqref{prob_init} by considering a set of candidate solutions, represented by particles. Then, particles move in the search space according to certain mathematical relations on the positions and velocities of the particles, until convergence toward a neighborhood of the global minimizer. Since PSO methods are gradient-free, the function $\mathcal{F}(x)$ in \eqref{prob_init} can be non-smooth.

The PSO has originally been formulated in such a way that each particle experiences a superposition of forces towards the local best position and towards the global best position found among all
particles up to that iteration \cite{kennedy,kennedy2}. There have been many extensions to the model; for example, in \cite{grassi2}, a description based on a system of stochastic differential equations has been introduced, as well as a formulation modeling memory of each particle. In another version of PSO, there is an inertia weight \cite{eberhart,shi} to bound the possible uncontrolled speed of the particles.
Additionally, PSO methods have gained prominence, drawing inspiration from CBO methods. In this variant of the PSO, the concept of a consensus point emerges, representing a weighted average of the positions of particles towards which they are attracted.

Several methods have been proposed in the literature for incorporating constraints in interacting particle consensus methods \cite{carrillo3, carrillo4, fornasier3, Fornasier4}. Here, we present a method based on exact penalization. This approach has already been used to study constrained problems with PSO methods in the version that uses the local best position and the global best position \cite{coello,Mezura-Montes2009}. However, in this manuscript we use the iterative strategy introduced for the CBO \cite{borghi2021constrained}, for which the penalty parameter is updated according to the violation of the constraint, and apply it to the version of the PSO method in which the consensus point appears.

The theoretical analysis of PSO methods has been developed on various scales \cite{grassi,grassi2}: the microscopic scale, the mesoscopic scale, and the macroscopic scale. The microscopic scale is based on the evolution of particle trajectories using systems of stochastic differential equations.
The intermediate scale, called mean-field, is the mesoscopic scale, where an equation governs the dynamics of the evolution of the probability density of the particles. These equations are a statistical description of the states of the particles.
Finally, at the macroscopic scale, the particle flow is considered to behave like a fluid. At this scale, the dynamics are governed by a system of hyperbolic conservation laws.

In this manuscript, we propose an approach called micro-macro decomposition; the idea is inspired by \cite{degond, lemoumehats, dimarcopareschi}. This decomposition is well-known in the literature in various contexts, for example, to derive the Navier-Stokes equations from the Boltzmann equation in rarefied gas dynamics \cite{Cercignani1988} or in the context of collisional kinetic equations that naturally incorporate exact space boundary conditions \cite{degond}. The decomposition is based on the decomposition of the probability density function into a contribution of the microscopic system and of the macroscopic system. The two contributions are weighted by a time-dependent function that is dynamically updated. This changes the behavior of the system to be more fluid-like or more particle-like. 

This manuscript is structured as follows. In Section \ref{sec:particle_swarm_optimization}, we present the penalization approach applied to constrained minimization problems. In Section \ref{sec:scale}, the microscopic model with inertia, its mean-field limit, and finally a second-order macroscopic system of the PSO methods are presented. In Section \ref{sec:micro_macro_decomposition}, we present the micro-macro decomposition. At last, in Section \ref{sec:numerical_results}, we present various numerical results obtained for constrained and unconstrained minimization problems, studying them at the microscopic and macroscopic scale and also using the mixed approach given by the micro-macro decomposition.

\section{Particle Swarm Optimization methods for constrained problems}
\label{sec:particle_swarm_optimization}
PSO methods are designed to solve unconstrained problems and we extend those to constrained problems of the following type:
\begin{equation*}
    \min_{x\in\mathbb{R}^d} \mathcal{F}(x) \quad \text{subject to } x\in \mathcal{M}
\end{equation*}
for $d\ge 1$, where $\mathcal{F}(x):\mathbb{R}^d\to \mathbb{R}$ is a given continuous, non-convex function, with an unique global minimizer in the feasible set $\mathcal{M}\subset\mathbb{R}^d$, which is closed, bounded, possibly non-convex, and empty. The set $\mathcal{M}$ is not required to be connected. In addition, we assume $\mathcal{M}$ has a boundary $\partial \mathcal{M}$ of zero Lebesgue-measure. In this paper we will focus on the study of the constrained case using a penalized PSO method, but other approaches are possible (see \cite{bonnans}). We consider exact penalization, since PSO does by itself not require differentiable objective functions.

The objective function is modified by a penalty term in such a way that points outside the set $\mathcal{M}$ are reached with low probability. The penalized problem to be solved is
\begin{equation}
\label{prob_pen}
    \min_{x\in\mathbb{R}^d}  \{\mathcal{F}_\beta(x):=\mathcal{F}(x)+\beta r(x) \},
\end{equation}
where $\beta \in \mathbb{R}^+$ is the strength of the penalization and $r$ is a penalty function, with the property:
\begin{equation*}
    r(x)=
    \begin{cases}
        =0\quad \text{if } x\in \mathcal{M},\\
        >0 \quad \text{if } x\not \in \mathcal{M}.\\
    \end{cases}
\end{equation*}
An example of $r(x)$ is the distance of the point $x\in\mathbb{R}^d$ from $\mathcal{M}$, e.g.,
\begin{equation}
\label{penalty_function_pso}
    r(x)=\text{min}_{y\in\mathcal{M}}||x-y||_p.
\end{equation}
The objective function $ \mathcal{F}_\beta$ coincides with the objective function $\mathcal{F}$ in the feasible set.

The goal is to obtain the same minimizer $x^*_{\mathcal{M}}$ for the penalized problem as for the constrained problem, for finite values of $\beta$. Penalization techniques that have this property are called exact penalization.
\begin{property}
\label{property1}
    A function $\mathcal{F}_\beta:\Omega \to \mathbb{R}$ is called an exact penalty function at a local minimum $x^*_{\mathcal{M}}$ of \eqref{prob_init}, if $x^*_{\mathcal{M}}$ is a local minimum of \eqref{prob_pen}\cite{bonnans,BERTSEKAS}.    
\end{property}
$F_\beta$ satisfies Property \ref{property1}, if $\beta$ is sufficiently large. If $r$ is non-differentiable, then $\beta$ is finite. We define this finite but a priori unknown value as $\Bar{\beta}$ 
\begin{equation*}
    \bar{\beta}:=\inf\Big\{\beta>0: \text{minimizers to problem } \eqref{prob_pen} \text{ are minimizers to problem } \eqref{prob_init}\Big\}.
\end{equation*} 
In the case of the PSO, as done for the CBO in \cite{borghi2021constrained}, we choose $\beta$ by an adaptive strategy \cite{herty}. The idea is to test whether the violation value $r(x)$ is less than or equal to a fixed tolerance. If the condition is satisfied, the feasibility requirement is increased. If the condition is not fulfilled, however, the penalty parameter $\beta$ is increased. The adaptive algorithm is explained in detail in Section \ref{constrained_alg}.\\
To summarize, the exact penalization was employed, thereby converting the initial problem \eqref{prob_init} into an unconstrained one \eqref{prob_pen}. This approach enables the employment of the penalized PSO method. The incorporation of the penalty function introduces a non-differentiability into the objective function. However, this does not pose a significant challenge for the PSO method, as it does not depend on derivatives, thereby negating the necessity for differentiability. An additional benefit is that by implementing an adaptive update to the parameter, a convergence to the solution of constrained problem \eqref{prob_init} is achieved. Furthermore, the method enables the resolution of constrained problems with unconnected constraints. In Section \ref{sec:numerical_results} various simulations are presented in which the feasible set $\mathcal{M}$ is not connected.

\section{The PSO method for constrained problems on microscopic and macroscopic scale}
\label{sec:scale}
\subsection{Microscopic scale}
In this section, we present the microscopic PSO method that describes the individual motion of particles, in the version proposed by \cite{grassi}.

The PSO is based on the classical second-order Newtonian dynamics of particle systems. At each time $n\geq 0$, we suppose to have $N$ particles, each with position $X^i_n\in\mathbb{R}^d$ and velocity $V^i_n\in\mathbb{R}^d$, $i=1,...,N$. We assign an initial position $X^i_0$ and an initial velocity $V^i_0$, then the evolution of positions and velocities is given by the canonical PSO equations 
\begin{equation}
\label{sist_particles}
\begin{cases}
    {X}^i_{n+1}={X}^i_{n}+\Delta  t  {V}^i_{n+1},\\
    {V}^i_{n+1}=\left(\frac{m}{m+\gamma\Delta t}\right) {V}^i_{n}+\left(\frac{\lambda \Delta t}{m+\gamma\Delta t}\right) ({X}^{\alpha,\beta}_n -{X}^i_n)+\left(\frac{\sigma \sqrt{\Delta t}}{m+\gamma\Delta t}\right)D({X}^{\alpha,\beta}_n - {X}^i_n) \theta^i_n,
\end{cases}
\quad \text{ for } i=1,...,N,
\end{equation}
with inertia weight $m \in (0,1]$, drift parameter $\lambda>0$, noise parameter $\sigma>0$, friction coefficient $\gamma:=(1-m)\geq0$, $\theta^i_n \sim \mathcal{N}(0,1)$ and $D({X}^{\alpha,\beta}_n - {X}^i_n)$ is
a $d$-dimensional matrix.
${X}^{\alpha,\beta}_n$ is the weighted average, commonly called consensus point, computed using the approach introduced for the CBO methods \cite{Pinnau} 
\begin{equation}
\label{sist_micro}
{X}^{\alpha,\beta}_n:=\frac{1}{N^{\alpha,\beta}_n}\sum_{i=1}^{N} {{X}}^i_n \omega^{\alpha,\beta}(X^i_n), \quad \text{with } N^{\alpha,\beta}_n:=\sum_{i=1}^N\omega^{\alpha,\beta}({X}^i_n).
\end{equation}
The weight function $\omega^{\alpha,\beta}$ is given by
\begin{equation*}
\omega^{\alpha,\beta}(X^{i}_n):=e^{{}-\alpha\mathcal{F_\beta}({{X}_n^i})}.
\end{equation*}
The choice of the Gibbs-type measure for $X^{\alpha,\beta}_n$ derived from the Laplace principle \cite{dembo}, which states that for any absolutely continuous distribution $\rho\in\mathcal{P}(\mathbb{R}^d)$ and $\mathcal{F}_\beta$ continuous, it holds
\begin{equation}
\label{laplace_princ}
    \lim_{\alpha\to\infty}\left(-\frac{1}{\alpha} \log \left(\int_{\mathbb{R}^d} e^{-\alpha\mathcal{F_\beta}(x)} d\rho(x)\right)\right)= \inf_{x\in supp(\rho)}\mathcal{F_\beta}.
\end{equation}
The proof of convergence of \eqref{laplace_princ} for functions with a unique global minimum is given e.g. in \cite{Pinnau}. By the Laplace principle \eqref{laplace_princ}, and for large values of $\alpha\gg 1$, the consensus point $X^{\alpha,\beta}_n$ is an approximation of the global minimizer of the objective function $\mathcal{F}_\beta$.


Various choices for the diffusion term $D({X}^{\alpha,\beta}_n - {X}^i_n)$ can be found in the literature. A first example is the isotropic exploration \cite{Pinnau}
\begin{equation*}
    D^{iso}({X}^{\alpha,\beta}_n - {X}^i_n)=||{X}^{\alpha,\beta}_n - {X}^i_n ||_2 I_d
\end{equation*}
where the norm used is the $L^2$ norm of vector in $\mathbb{R}^d$ and $I_d\in \mathbb{R}^{d\times d}$ is the identity matrix. Another possibility is the anisotropic exploration \cite{Carrillo2}
\begin{equation}
\label{anisotropic}
        D^{aniso}({X}^{\alpha,\beta}_n - {X}^i_n)=\text{diag}{(({X}^{\alpha,\beta}_n - {X}^i_n)_1,...,({X}^{\alpha,\beta}_n - {X}^i_n)_d)}
 \end{equation}
where $\text{diag}:\mathbb{R}^d \to \mathbb{R}^{d\times d}$ is the operator that maps a vector onto a diagonal matrix. For high-dimensional problems, it has been shown to be more efficient to use anisotropic exploration \cite{Carrillo2,fornasier2} and we use $D({X}^{\alpha,\beta}_n - {X}^i_n)=D^{aniso}({X}^{\alpha,\beta}_n - {X}^i_n)$.

The change of the velocity in \eqref{sist_particles} is subject to three forces. The first term is the friction term. It depends only on the speed of the $i$-th particle in the previous time step. The second term is the exploitation term, which models the drift of the $i$-th particle towards the consensus point $X^{\alpha,\beta}_n$. The contribution of exploitation to acceleration depends on the distance of the particle from $X^{\alpha,\beta}_n$. The third term is the exploration term and it is the only random contribution.

The system \eqref{sist_particles} is an Euler-Maruyama time discretization of the system of stochastic differential equations (SDEs) in It\^{o} form \cite{higham}
\begin{equation*}
\begin{cases}
    dX_t^i=V_t^i dt,\\
    dV_t^i=-\frac{\gamma}{m} V_t^i dt+\frac{\lambda}{m}(X^{\alpha,\beta}_{t} -X_t^i)dt+\frac{\sigma}{m} D(X^{\alpha,\beta}_t-X_t^i)dB_t^i,
\end{cases}
\quad \text{ for } i=1,...,N,
\end{equation*}
where $X^{\alpha,\beta}_{t} $ is defined as \eqref{sist_micro} and $\{(B_t^i)_{t\geq0}\}_{i=1}^N$ are $N$ independent $d$-dimensional standard Brownian motions in $\mathbb{R}^d$. 
 
\subsection{Mean-field limit and macroscopic scale}
We derive a mean-field formulation and the macroscopic system associated with the system \eqref{sist_particles}. On the mesoscopic scale, a probabilistic approach is considered through a distribution function $f=f(t,x,v)\in\mathcal{P}(\mathbb{R}^+\times\mathbb{R}^d\times\mathbb{R}^d)$ that represents the density of particles in the phase-space $(x,v)$. The behavior of the particles is predicted through the study of statistical quantities of $f$. For the derivation, refer to \cite{Huang_2022}.\\
The mean-field equation of the system \eqref{sist_particles} is
\begin{equation}
\label{mean_field}
    \partial_t{f}+v\cdot \nabla_x{f}=\nabla_v \cdot \left(\frac{\gamma}{m}vf+\frac{\lambda}{m}(x-X_{\alpha,\beta}(\rho))f+\frac{\sigma^2}{2m^2}D(x-X_{\alpha,\beta}(\rho))^2\nabla_v \cdot {f}\right)
\end{equation}
where
\begin{equation}
\label{eqn:consensus_point_macro}
    X_{\alpha,\beta}(\rho):=\frac{\int_{\mathbb{R}^d}x \omega_{\alpha,\beta} (x) \rho(t,x)dx}{\int_{\mathbb{R}^d} \omega_{\alpha,\beta} (x) \rho(t,x)dx} \quad \text{and} \quad \rho=\rho(t,x)=\int_{\mathbb{R}^d} f(t,x,v)dv.
\end{equation}
The mean-field equation can be written in compact form as
\begin{equation*}
    \partial_t{f}+v\cdot \nabla_x{f}+Q[f]=0,
\end{equation*}
where 
\begin{equation*}
Q[f]:=-\nabla_v \cdot \left(\frac{\gamma}{m}vf+\frac{\lambda}{m}(x-X_{\alpha,\beta}(\rho))f+\frac{\sigma^2}{2m^2}D(x-X_{\alpha,\beta}(\rho))^2\nabla_v \cdot{f}\right)
\end{equation*}
is the collisional operator. It acts as a non-local operator on $f$ only in velocity and models how binary collisions between particles affect the evolution of the distribution function $f$, accounting for both the gain and loss of particles in a given velocity state due to these interactions \cite{Cercignani1988}. We observe that the collisional operator $Q$ satisfies the local conservation property for $\phi(v):=\left(\begin{matrix}
    & 1 &\\
    & v & 
\end{matrix}\right)$, i.e.,
\begin{equation*}
    \langle \phi(v),Q(g)\rangle=0 \quad \quad \forall g\in\mathcal{P}(\mathbb{R}^+\times\mathbb{R}^d\times\mathbb{R}^d).
\end{equation*}
The components of the vector $\phi$ are the collisional invariants, namely the locally conserved quantities in $(t,x)$. Therefore we define the vector $U=U(t,x)$ as
\begin{equation*}
    U(t,x)=\left(\begin{matrix}
    & \rho &\\
    & \rho u& 
\end{matrix}\right)(t,x):=\int_{\mathbb{R}^d} f(t,x,v)
 \left(   \begin{matrix}
    & 1 &\\
    & v & 
\end{matrix}\right)dv.
\end{equation*}
To obtain equations for $U$ of the macroscopic system, we integrate the mean-field equation \eqref{mean_field} against the collisional invariants $\phi(v)$ and obtain\\
\begin{equation}
\label{eqn: system_no_closure}
\partial_t \left(
\begin{matrix}
    & \rho &\\
    & \rho u& 
\end{matrix}
\right)+\nabla_x \cdot \left(
\begin{matrix}
    & \rho u&\\
    &\int v^2 f(t,x,v) dv& 
\end{matrix}
\right)+\left(
\begin{matrix}
    & 0 &\\
    & \frac{\gamma}{m}(\rho u)+\frac{\lambda}{m}(x-X_{\alpha,\beta}(\rho))\rho& 
\end{matrix}
\right)=0.
\end{equation}
The initial data $U_0(x)=U(0,x)$ is given by
\begin{equation*}
U_0(x)=\left(\begin{matrix}
    & \rho_0 &\\
    & (\rho u)_0& 
\end{matrix}\right)(x)=\int_{\mathbb{R}^d} f(0,x,v)
 \phi(v)dv.
\end{equation*}\\
In its compact form, system \eqref{eqn: system_no_closure} is written as
\begin{equation*}
    \partial_t U+\nabla_x \cdot F(U) +S(U)=0,
\end{equation*}
where \begin{equation}
\label{eqn: def_U}
    U:=\left(
\begin{matrix}
    & \rho &\\
    & \rho u& 
\end{matrix}
\right)\text{ is the vector of conserved variables, }
\end{equation}
\begin{equation*}
F(U):=\left(
\begin{matrix}
    & \rho u&\\
    &\int v^2 f(t,x,v) dv& 
\end{matrix}
\right) \text{ is the flux vector,}
\end{equation*}
\begin{equation}
\text{ and } S(U):=\left(
\begin{matrix}
\label{eqn: def_source_term}
    & 0 &\\
    & \frac{\gamma}{m}\rho u+\frac{\lambda}{m}(x-X_{\alpha,\beta}(\rho))\rho& 
\end{matrix}
\right) \text{ is the source term},
\end{equation}
with $X_{\alpha,\beta}(\rho)$ defined as in \eqref{eqn:consensus_point_macro}.
The derivation of the macroscopic system is conducted by first starting from the mean-field equation and then exploiting the properties of fluid dynamics, as referenced in \cite{bardos1, Cercignani1988}. \\
The system \eqref{eqn: system_no_closure} depends on the distribution $f$. To eliminate this dependence, a closure relation can be imposed by selecting a specific distribution, $f^*=f^*(t,x,v)$, which ensures that the conservation laws are hyperbolic or strictly hyperbolic. By expressing the system \eqref{eqn: system_no_closure} in quasi-linear form, it is possible to determine the properties that the function $f^*$ must satisfy. \\
In order to derive the quasi-linear form, it is necessary to express the spatial derivative of the flows $\nabla_x \cdot F(U)$, as a function of the derivatives of $U$. To accomplish this, the chain rule is applied component by component
\begin{equation*}
    \frac{\partial F_i(U)}{\partial x_i}=\frac{\partial F_i(U)}{\partial U} \partial_{x_i} U=A_i(U)\partial_{x_i} U, \quad \quad i=1,...,d.
\end{equation*}
For the system \eqref{eqn: system_no_closure}, we obtain the following:
\begin{equation}
\label{eqn: def_Ai}
    A_i(U):=\frac{\partial}{\partial U} F_i(U)=\frac{\partial}{\partial U}\left(
\begin{matrix}
    & \rho u&\\
    &\int v^2 f(t,x,v) dv& 
\end{matrix}
\right)_i= \left(
\begin{matrix}
    & 0& e_i^T&\\
    & \partial_\rho \int v^2 f^* dv& \partial_{\rho u} \int v^2 f^* dv&  
\end{matrix}
\right).
\end{equation}
Accordingly, the system \eqref{eqn: system_no_closure} is represented in quasi-linear form as
\begin{equation}
\label{eqn: system_no_closure_quasilin}
    \partial_t U+\sum_{i=1}^d A_i(U) \partial_{x_{i} }U+S(U)=0
\end{equation}
where $U$, $S(U)$ and $A_i(U)$ are defined in \eqref{eqn: def_U}, \eqref{eqn: def_source_term} and \eqref{eqn: def_Ai}, respectively.\\
In order to obtain a well-posed initial value problem, we introduce a basic algebraic property.
\begin{definition}
    A system of balance laws
    \begin{equation*}
       \partial_t U+\sum_{i=1}^d A_i(U) \partial_{x_{i} }U+S(U)=0
    \end{equation*}
    in $\mathbb{R}^n$, with initial data $U(x,0)=U_0(x)$, is strictly hyperbolic if, for every $U$, the composite Jacobian matrix in the direction $\xi$
    \begin{equation*}
        A(\xi)=\sum_{i=1}^d A_i(U) \xi_i= \sum_{i=1}^{d}\frac{\partial F_i}{\partial U}\xi_i
    \end{equation*}
    has $n$ real and distinct eigenvalues: $\lambda_1(U)<...<\lambda_n(U)$.
\end{definition}
From now on, we focus on problems in dimension 1. It is possible to derive a sufficient condition for the system \eqref{eqn: system_no_closure_quasilin} to be strictly hyperbolic.
\begin{lemma}
In dimension 1, the system is strictly hyperbolic if the following condition is satisfied
\begin{equation}
\label{eqn:cond_hyperbolicity}
    \left(\partial_{\rho u} \int v^2 f^* dv\right)^2>4\partial_\rho \int v^2 f^* dv.
\end{equation}
\end{lemma}
\begin{proof}
In order to identify a sufficient condition for strict hyperbolicity, which $f^*$ must satisfy, we compute the eigenvalues of the matrix $A(U)$. In dimension 1 the system in quasi-linear form is 
\begin{equation}
\label{eqn: quasilinear1d}
    \partial_t U+A(U)\partial_x U+S(U)=0, \quad \text{with} \quad A(U)=\left(
\begin{matrix}
    & 0& 1&\\
    & \partial_\rho \int v^2 f^* dv& \partial_{\rho u} \int v^2 f^* dv &  
    \end{matrix} \right).
\end{equation}
We proceed by calculating the eigenvalues of $A(U)$
\begin{equation*}
\begin{aligned}
    det (A(U)-\lambda Id)&=0,\\
    det \left(
\begin{matrix}
    & -\lambda& 1&\\
    & \partial_\rho \int v^2 f^* dv& \partial_{\rho u} \int v^2 f^* dv -\lambda&  
    \end{matrix} \right) &=0,\\
    \lambda^2-\partial_{\rho u} \int v^2 f^* dv \lambda-\partial_\rho \int v^2 f^* dv &=0.\\
\end{aligned} 
\end{equation*} 
To obtain real and distinct eigenvalues, it is necessary that the discriminant of the last quadratic equation be strictly positive, i.e.
\begin{equation*}
    \left(\partial_{\rho u} \int v^2 f^* dv\right)^2-4\partial_\rho \int v^2 f^* dv>0.
\end{equation*}
It thus follows that in order for the system \eqref{eqn: quasilinear1d} to be strictly hyperbolic and therefore well-posed, the distribution $f^*$ selected for closure must satisfy the following condition:
\begin{equation*}
    \left(\partial_{\rho u} \int v^2 f^* dv\right)^2>4\partial_\rho \int v^2 f^* dv.
\end{equation*}
\end{proof}

In dimension 1, a potential selection for $f^*$ that satisfies \eqref{eqn:cond_hyperbolicity} is the Maxwellian distribution with parameter $T\in \mathbb{R}^+$
\begin{equation}
\label{eqn: macrosyst_closure}
f^*(t,x,v)=\frac{\rho}{\sqrt{2\pi T^2}} \exp{\left(-\frac{1}{2}\frac{(v-u(t,x))^2}{T^2}\right)},
\end{equation}
where $u(t,x):=\frac{\rho u(t,x)}{\rho(t,x)}$ indicates the average value of the velocity along the direction $v$, at time $t$ in $x$, and $T^2$ is the variance of the velocity. $T$ is chosen in $\mathbb{R}^+$ because it represents the standard deviation.\\
Using the expected value and variance properties of the normal distribution, we obtain $\int v^2 f^* dv = \rho u^2 +\rho T^2 $.
In this case the system of differential equations for $U$ is
\begin{equation}
\label{eqn: macroscopic_system_2d}
\partial_t \left(
\begin{matrix}
    & \rho &\\
    & \rho u& 
\end{matrix}
\right)+\partial_x \left(
\begin{matrix}
    & \rho u&\\
    & \rho u^2 +\rho T^2 & 
\end{matrix}
\right)+\left(
\begin{matrix}
    & 0 &\\
    & \frac{\gamma}{m}(\rho u)+\frac{\lambda}{m}(x-X_{\alpha,\beta}(\rho))\rho& 
\end{matrix}
\right)=0
\end{equation}
and with initial data $U_0$
\begin{equation*}
U_0(x)=\left(\begin{matrix}
    & \rho_0 &\\
    & (\rho u)_0& 
\end{matrix}\right)(x)=\int_{\mathbb{R}^d} f(0,x,v)
 \phi(v)dv.
\end{equation*}\\
The associated quasi-linear system can be expressed as follows:
\begin{equation}
\label{eqn: macroscopic_system_2d_quasilin}
    \partial_t U+A(U)\partial_x U+S(U)=0
\end{equation}
with $U=\left(
\begin{matrix}
    & \rho &\\
    & \rho u& 
\end{matrix}
\right)$, $A(U)=\left(
\begin{matrix}
    & 0& 1&\\
    & T^2-u^2& 2u& 
\end{matrix}
\right) $ and $S(U)=\left(
\begin{matrix}
    & 0 &\\
    & \frac{\gamma}{m}(\rho u)+\frac{\lambda}{m}(x-X_{\alpha,\beta}(\rho))\rho& 
\end{matrix}
\right)$. 
The eigenvalues of the matrix $A(U)$ are given by
\begin{equation*}
    \lambda_1=u+|T| \quad \text{and} \quad \lambda_2=u-|T|.
\end{equation*}
The system \eqref{eqn: macroscopic_system_2d_quasilin} is strictly hyperbolic for $T\not=0$, which implies that the system \eqref{eqn: macroscopic_system_2d} is well posed \cite{kreiss}.

\begin{remark}
    It should be noted that the monokinetic closure $f^*=\rho \delta(v-u)$ does not guarantee the strict hyperbolicity of the system, as it does not satisfy condition \eqref{eqn:cond_hyperbolicity}. Furthermore, it does not guarantee hyperbolicity, as the matrix $A(U)=\left(
\begin{matrix}
    & 0 & 1&\\
    & -(\rho u)^2/\rho^2& 2(\rho u)/\rho& 
\end{matrix}
\right)$, obtained in the quasi-linear form, is not diagonalizable.
\end{remark}

\begin{remark}
It is also possible to derive a third-order macroscopic system. In this case we will have that the vector of conserved quantities is defined as $\phi(v):=\left(\begin{matrix}
    & 1 &\\
    & v & \\
    & v^2&
\end{matrix}\right)$ and the vector $U$ is defined as
\begin{equation*}
    U(t,x)=\left(\begin{matrix}
    & \rho &\\
    & \rho u&\\
    & \rho T^2&
\end{matrix}\right)(t,x):=\int_{\mathbb{R}^d} f(t,x,v)
 \left(   \begin{matrix}
    & 1 &\\
    & v &\\
    & v^2&
\end{matrix}\right)dv.
\end{equation*}
Again, the Maxwellian distribution \eqref{eqn: macrosyst_closure} can be chosen as the closure.
Given $U_0$ initial data, we obtain the third-order macroscopic system in quasi-linear form as
\begin{equation*}
    \partial_t U+A(U)\partial_x U+S(U)=0,
\end{equation*}
	\begin{equation*}
	\text{with} \quad U=\left(
	\begin{matrix}
		& \rho &\\
		& \rho u& \\
		& \rho T^2&
	\end{matrix}
	\right), \quad  A(U)=\left(
	\begin{matrix}
		& 0& 1& 0&\\
		& -u^2& 2u& 1&\\
		&-3uT^2& 3T^2& u&
	\end{matrix}
	\right)  \quad \text{and} 
	\end{equation*}
	\begin{equation*}
		S(U)=\left(
		\begin{matrix}
			& 0 &\\
			& \frac{\gamma}{m}(\rho u)+\frac{\lambda}{m}(x-X_{\alpha,\beta}(\rho))\rho& \\
			&2\frac{\gamma}{m}\rho T^2
		\end{matrix}
		\right). 
	\end{equation*}
The system is strictly hyperbolic since the eigenvalues are real and distinct
\begin{equation*}
    \lambda_1=u, \quad \lambda_2=u+\sqrt{3}T \quad \text{and} \quad \lambda_3=u-\sqrt{3}T.
\end{equation*}
\end{remark}
In conclusion, two formulations were proposed to solve the problems presented in \eqref{prob_init} and in \eqref{prob_pen}. The initial approach was the microscopic system \eqref{sist_particles}, which is based on the dynamics of individual particles. The evolution of the particles is given by the evolution of their positions and velocities at each time $n$. The velocity is updated based on the previous velocity of the particle and the distance of the particle from the consensus point, which is calculated at each time instant. The convergence of this system is governed by the relationship $X_{\alpha,\beta}\to argmin \mathcal{F}_{\beta}$. 

The second approach was to solve the minimization problem using the macroscopic system \eqref{eqn: macroscopic_system_2d}, which is based on the evolution of $U$. We also derived the quasi-linear formulation of these conservation laws. In addition, we proved a criterion that the closure distributions must meet in dimension 1 to guarantee the strict hyperbolicity of the system. In the simulations in Section \ref{sec:numerical_results}, tests obtained by implementing the conservation laws in dimension 1 and with the closure relation obtained using the Maxwellian distribution \eqref{eqn: macrosyst_closure} are presented. In this case, convergence is given by the concentration of the mass of the system around the minimizer $x^*_{\mathcal{M}}$, i.e. $supp(\rho)\to argmin \mathcal{F}_{\beta}$.

\section{Micro-macro decomposition} \label{sec:micro_macro_decomposition}
We propose an approach based on a formal micro-macro decomposition of the kinetic probability density to exploit the scales. The approach develops the idea of studying and comparing the behavior of the PSO method applied to the microscopic scale and the macroscopic scale, which has been presented in Section \ref{sec:scale}. Specifically, the decomposition is employed to calculate the two dynamics (microscopic and macroscopic) and couple them, with the objective of efficiently moving mass between the two scales to achieve method convergence. \\The micro-macro decomposition typically leads to more efficient computational methods by reducing the number of particles to be simulated. Additionally, merging the two dynamics leverages the fact that the macroscopic model is defined across the entire spatial domain. This property enables the computation of consensus points that are a good approximation of the minimizer $x^*_{\mathcal{M}}$, even for small times.

The idea of this approach is inspired by \cite{degond, lemoumehats, dimarcopareschi}, where micro-macro decomposition is used to decompose the spatial domain to match fluid and kinetic models. Here, we use the micro-macro approach directly on the distribution function $f=f(t,x,v)\in\mathcal{P}(\mathbb{R}^+\times\mathbb{R}^d\times\mathbb{R}^d)$, introduced in \eqref{mean_field}.\\
The starting point is to formally split the distribution function $f$ in \eqref{mean_field} in a microscopic contribution $\hat{f}=\hat{f}(t,x,v)\in\mathcal{P}(\mathbb{R}^+\times\mathbb{R}^d\times\mathbb{R}^d)$ and a macroscopic contribution, specifically, we write using the closure of the macroscopic system $f^*=f^*(t,x,v)\in\mathcal{P}(\mathbb{R}^+\times\mathbb{R}^d\times\mathbb{R}^d)$ \eqref{eqn: macrosyst_closure}
\begin{equation*}
f(t,x,v)=\zeta(t,x) \hat{f}(t,x,v)+(1-\zeta(t,x))f^*(t,x,v), \quad\quad \zeta(t,x)\in [0,1],
\end{equation*}
where $(x,v)\in (\mathbb{R}^d\times\mathbb{R}^d)$. To compute the evolution of the distribution, we do not solve the mean-field equation but instead use particle dynamics and the macroscopic system. At the microscopic level the evolving quantities are the positions and velocities of the particles, so the microscopic contribution is given by the empirical distribution associated with the system \eqref{sist_particles}
\begin{equation*}
    \hat{f}(t,x,v)=\frac{1}{N}\sum_{i=1}^N \delta (x-X^i_n) \delta (v-V^i_n).
\end{equation*}
Instead, the macroscopic contribution is made by exploiting the closure relation 
\begin{equation*}
f^*(t,x,v)=\frac{\rho (t,x)}{\sqrt{2\pi T^2}} \exp{\left(-\frac{1}{2}\frac{(v-u(t,x))^2}{T^2}\right)},
\end{equation*}
which requires to propagate $(\rho, \rho u)^T$ accordingly to \eqref{eqn: macroscopic_system_2d}. The weight of the total distribution $f$ allocated to the microscopic and macroscopic scales is given by the function $\zeta(t,x)$. Therefore, the role of $\zeta(t,x)$ is to establish a connection between the two dynamics. The coupling is determined by the movement of mass shift between the microscopic and macroscopic scales. There are several possible choices for this function. 
The approach employed in this study involves the partitioning of the spatial domain into cells. The $j$-th cell of the spatial grid is denoted by $I_j$. The value of the function $\zeta$ is given by

\begin{equation}
\label{def: zeta}
\zeta(t,x):=\frac{\sum_{j=1}^K w(t,I_j) |u(t,I_j)-\sum_{i=1}^N\mathbbm{1}_{I_j}(X_t^i)\cdot V_t^i|}{\sum_{j=1}^K w(t,I_j) \max_{i=1,...,N}|u(t,I_j)-\mathbbm{1}_{I_j}(X_t^i)\cdot V_t^i|}, \quad \forall x \in I_j,
\end{equation} 
where $K$ is the number of cells, $N$ is the number of particles, $u(t,x)$ is defined as 
\begin{equation*}
    u(t,I_j):=\int{ I_j}\frac{\rho u(t,x)}{\rho(t,x)}dx
    \end{equation*}
and $V_t^i$ are the velocities of the particles. The weight $w(t,I_j)$ is defined as the relative density
\begin{equation*}
w(t,I_j):=\int_{ I_j}\frac{\rho^{m}(t,x)}{\rho^{M}(t,x)+\rho^{m}(t,x)}dx,
\end{equation*}
where
\begin{equation}
\label{def: rho_micro}
\rho^{m}(t,x)=\int_{\mathbb{R}^d} \hat{f}(t,x,v)dv \quad \text{ and }\quad \rho^{M}(t,x)=\int_{\mathbb{R}^d} f^*(t,x,v)dv, \quad \text{respectively.}
\end{equation}
The definition of $\zeta$ is motivated by the following considerations. The difference $|u(t,I_j)-\sum_{i=1}^N\mathbbm{1}_{I_j}(X_t^i)\cdot V_t^i|$  measures the distance between the macroscopic velocity defined by $u$ \eqref{eqn: macroscopic_system_2d} and the average velocity of the particles in the same cell \eqref{sist_particles}. This quantity is significant since if microscopic velocities tend towards the chosen velocity of the macroscopic system, they become less relevant. Hence, we move mass to macroscopic scale. The distance is weighted by the local mass $w(t,I_j)$. Therefore, the effect of distance is determined by the microscopic mass in the cell $I_j$, and its contribution becomes more significant as the local mass increases. The denominator normalizes the quantity $\zeta$. \\
The function $\zeta$ balances the total density $\rho(t,x)=\int_{\mathbb{R}^d} f(t,x,v)dv$ of the system. This balance is based on the comparison of the output of densities obtained at each scale at each time instant. If the microscopic system gives results that are far from those obtained by the macroscopic system, $\zeta$ will take large values. Then, the mass is shifted to the microscopic scale that is expected to be more accurate. In the opposite case, where the particle velocities have values close to those of the macroscopic velocity, the microscopic system is expected to be closer to the macroscopic system. Then, we expect a higher density at the macroscopic scale. \\
The detailed algorithm is given in Section \ref{constrained_alg}.

\section{Numerical algorithm and results} \label{sec:numerical_results}
In this section, several simulations are presented. The numerical results obtained using the algorithm introduced at the microscopic and macroscopic scales are reported, both in the classical and in the constrained case. In addition, results are given for the micro-macro decomposition. 

\subsection{Constrained micro-macro PSO method} \label{constrained_alg}
In Section \ref{sec:particle_swarm_optimization} we introduced the concept of adaptive exact penalty in order to solve constrained minimization problems. This strategy succeeds in handling even non-differentiable functions. The idea is to minimize a modified objective function $\mathcal{F}_\beta$ over the entire space $\mathbb{R}^d$, i.e.
\begin{equation*}
    \min_{x\in\mathbb{R}^d} \mathcal{F}_\beta(x):=\mathcal{F}(x)+\beta r(x),
\end{equation*}
where $\beta \in \mathbb{R}^+$ and $r$ is the penalty function, in our simulations is given by
\begin{equation*}
    r(x)=\text{min}_{y\in\mathcal{M}}||x-y||_p.
\end{equation*}
In this study, an iterative algorithm in $\beta$ is considered, and the dependence on the interaction $n$ is denoted by $\beta_n$. The evolution of the parameter $\beta_n$ is based on the violation of the constraint $r$ \eqref{penalty_function_pso} by the microscopic/macroscopic system. This violation can be measured on the microscopic scale by 
\begin{equation}
\label{R_micro}
    \mathcal{R}^{m}(X_n):=\frac{ \sum_{i=1}^{N} r(X^i_n) \exp(-\alpha \mathcal{F}_\beta(X^i_n))}{\sum_{i=1}^{N} \exp(-\alpha \mathcal{F}_\beta(X^i_n))}
\end{equation}
and on the macroscopic scale by   
\begin{equation*}
\label{R_macro}
    \mathcal{R}^{M}(x):= \frac{\int r(x) \exp(-\alpha \mathcal{F}_\beta(x))\rho(x) dx}{\int \exp(-\alpha \mathcal{F}_\beta(x))\rho(x) dx}, \text{ respectively.}
\end{equation*}
The violation is measured with a weighted expectation of $r$ based on the Boltzmann-Gibbs distribution, the same one used for the calculation of the consensus point \eqref{sist_micro} (respectively \eqref{eqn:consensus_point_macro}). This choice results in an adaptive mechanism that is more accurate than that obtained with other feasibility conditions, as shown in \cite{borghi2021constrained}.\\
In addition, a tolerance value, denoted by $\ka_n$, is introduced and subject to temporal variation. This value is used as a threshold for the constraint violations. In each temporal interval, it is necessary to ascertain whether the following feasibility condition is satisfied:
\begin{equation}
\label{eqn: feasibility_check}
   \mathcal{R}^{m}\leq \frac{1}{\sqrt{\ka_n}} \quad \text{ or respectively  } \quad
\mathcal{R}^{M}\leq \frac{1}{\sqrt{\ka_n}} .
\end{equation}
The formulation of the adaptive strategy is founded upon the evaluation of the validity of the inequality \eqref{eqn: feasibility_check}, which is tested at each time step, resulting in updated values of $\beta_n$ and $\kappa_n$. If the condition \eqref{eqn: feasibility_check} holds, it means that the violation of the constraint is contained within the imposed value $1/\sqrt{\ka_n}$. Consequently, the variation of the parameter $\beta_n$, which governs the penalty contribution, is avoided
\begin{equation*}
        \beta_{n+1}=\beta_n.
    \end{equation*}
However, the necessity arises to update the value of $\kappa_n$ in order to test the inequality at the next time step under more restrictive conditions. The value $\kappa_{n+1}$ is updated in the following way
    \begin{equation*}
    \ka_{n+1}=\eta_{\ka}\ka_n, \quad \text{ with $\eta_\ka>1$. }
    \end{equation*}
On the other hand, if the condition \eqref{eqn: feasibility_check} is not satisfied, it means that the contribution of the violation $\mathcal{R}^M$ exceeds the established tolerance $\kappa_n$. It is therefore necessary to increase the parameter $\beta_{n+1}$ in order to impose additional particles (or mass density, respectively) in the constraint 
\begin{equation*}
        \beta_{n+1}=\eta_\beta\beta_n, \quad \text{ with $\eta_\beta>1$. }
    \end{equation*}
The result is that the contribution of $\mathcal{R}^M$ is decreased. Consequently, the value of $\kappa_{n+1}$ is reduced in order to test the inequality with a lower tolerance at the subsequent time. Therefore,
\begin{equation*}
        \ka_{n+1}=\min\{\ka_{n}/\eta_\ka,\ka_0\}, \quad \text{ with $\eta_\ka>1$. }
    \end{equation*}

The strategy is presented in Algorithm \ref{micromacro_pseudocode}.

\begin{algorithm}
\caption{Algorithm for the micro-macro decomposition in the constrained case}
\begin{algorithmic} [1]
\label{micromacro_pseudocode}
\REQUIRE \text{initial data }$N_t, \alpha, \beta, X_0$, $\rho_0$, $ \zeta_0$, $\mu_0$
\REQUIRE \text{parameters of the penalized microscopic algorithm }$\kappa_0^m, \eta_\kappa^m, \beta_0^m, \eta_\beta^m$
\REQUIRE \text{parameters of the penalized macroscopic algorithm }$\kappa_0^M, \eta_\kappa^M, \beta_0^M, \eta_\beta^M$
\STATE \text{compute the density $\rho^{m}_0$ associated to the microscopic scale \eqref{def: rho_micro}}
\FOR{$n=1,...,N_t$}
\STATE \% penalized microscopic algorithm 
\STATE apply the Euler-Maruyama scheme  \eqref{sist_particles}
\STATE $X_n^{\alpha,\beta} = \frac{1}{N^{\alpha,\beta}_n}\sum_{i=1}^{N} {{X}}^i_n \omega^{\alpha,\beta}(X^i_n)$ 
\IF{\text{microscopic feasibility condition} \eqref{eqn: feasibility_check} \text{ holds}}
\STATE $\kappa^{m }_{n+1}=\eta^{m }_\kappa \kappa^{m }_n$
\STATE $\beta^{m }_{n+1} = \beta^{m }_n$
\ELSE
\STATE  $\ka^{m}_{n+1}=\min\{\ka^{m}_{n}/\eta^{m}_\ka,\ka^{m}_0\}$
\STATE $\beta^{m}_{n+1} =\eta^{m}_\beta \beta^{m}_n$
\ENDIF
\STATE \% penalized macroscopic algorithm 
\STATE apply the Lax-Friedrichs scheme \eqref{eqn: LxF} for the system \eqref{eqn: macroscopic_system_2d}
\STATE $X_{\alpha,\beta}(\rho) = \frac{\int_{\mathbb{R}^d}x \omega_{\alpha,\beta} (x) \rho(t,x)dx}{\int_{\mathbb{R}^d} \omega_{\alpha,\beta} (x) \rho(t,x)dx} $,
\IF{\text{macroscopic feasibility condition} \eqref{eqn: feasibility_check}\text{ holds}}
\STATE $\kappa^{M }_{n+1} =\eta^{M }_\kappa \kappa^{M }_n$
\STATE $\beta^{M }_{n+1} =\beta^{M }_n$
\ELSE
\STATE  $\ka^{M}_{n+1}=\min\{\ka^{M}_{n}/\eta^{M}_\ka,\ka^{M}_0\}$
\STATE $\beta^{M}_{n+1} = \eta^{M}_\beta \beta^{M}_n$
\ENDIF
\STATE \% micro-macro algorithm 
\STATE compute $\zeta_n$ according to \eqref{eqn: zeta_numeric}
\STATE $\mu_n=\zeta_n\cdot \mu_0$ $\quad$ \% update of the mass of particles ($\mu_0$ is the initial mass of particles) 
\STATE compute the density $\rho^{m}_n$ associated to the microscopic scale \eqref{def: rho_micro}
\IF {$\mu_n\leq\mu_{n-1}$} 
\STATE $\rho^{M}_n = \rho^{M}_n+(\rho^{m}_n-\rho^{m}_{n-1})$
\ELSE
\STATE $\rho^{M}_n =\rho^{M}_n-(\rho^{m}_n-\rho^{m}_{n-1})$
\ENDIF
\ENDFOR

\end{algorithmic}

\end{algorithm}

\begin{remark}
The iterative strategy previously outlined has been implemented and its convergence is demonstrated for another meta-heuristic method for solving minimization problems, namely the Consensus-Based Optimization method, as detailed in \cite{borghi2021constrained}.
\end{remark}
The penalized adaptive strategy presented for solving constrained problems can be used, as already mentioned, on both the microscopic and macroscopic scales. In the former case, the positions and velocities of the particles are calculated using the Euler-Maruyama scheme \eqref{sist_particles}. Instead, on the macroscopic scale, the quantity that evolves is $U$. To compute it, we discretize the system \eqref{eqn: macroscopic_system_2d} using the Lax-Friedrichs scheme 
\begin{equation}
\label{eqn: LxF}
U^{\ell}_{n+1}=\frac{1}{2}(U^{\ell+1}_n+U^{\ell-1}_n)-\frac{\Delta t}{2\Delta x}(F(U^{\ell+1}_n)-F(U^{\ell-1}_n))+\Delta t S(U_n^{\ell}), 
\end{equation}
with $\ell=1,...,L$ are the grid indices on which the scheme is applied, where the total number of nodes is $L$, $U^{\ell}_n:=\left(
\begin{matrix}
    & \rho^{\ell}_n &\\
    & \rho u^{\ell}_n& 
\end{matrix}
\right)$, $F(U^{\ell}_n):=\left(
\begin{matrix}
    & \rho u^{\ell}_n&\\
    & (\rho u^2 +\rho T^2)^{\ell}_n & 
\end{matrix}
\right)$, $S(U^{\ell}_n):=\left(
\begin{matrix}
    & 0 &\\
    & \frac{\gamma}{m}\rho u^{\ell}_n+\frac{\lambda}{m}(x^{\ell}-X_{\alpha,\beta}(\rho^{\ell}_n))\rho^{\ell}_{n+1}& 
\end{matrix}
\right) $ and
$x_{\ell}$ are the points of the spatial grid defined on the domain of the objective function $\mathcal{F}_{\beta}$.

We proceed to present the algorithm for micro-macro decomposition. Consider the local densities defined over a domain. Then, introduce a spatial grid over the domain, dividing it into cells. The $j$-th cell of the spatial grid is denoted by $I_j$.
The value of the function $\zeta$ is computed of each cell in the grid, so that
\begin{equation}
\label{eqn: zeta_numeric}
\zeta(t,I_j):=\frac{\sum_{j=1}^K\sum_{i=1}^N w(t,I_j) |u(t,I_j)-\mathbbm{1}_{I_j}(X_t^i)\cdot V_t^i|}{\sum_{j=1}^K w(t,I_j) \max_{i=1,...,N}|u(t,I_j)-\mathbbm{1}_{I_j}(X_t^i)\cdot V_t^i|}.
\end{equation} 
where $ u(t,I_j):=\int{ I_j}\frac{\rho u(t,x)}{\rho (t,x)}dx$ and $w(t,I_j):=\int_{ I_j}\frac{\rho^{m}(t,x)}{\rho^{M}(t,x)+\rho^{m}(t,x)}dx$. To compute the quantities $u$ and $w$, we used a first-order quadrature rule.

The basic steps of micro-macro decomposition, introduced for the PSO method, are as follows:
\begin{itemize}
\item fix the initial mass of the particles $\rho(0,x)=\rho^{m}(0,x)+\rho^{M}(0,x)$;
\item at each time step compute the new value of the particle mass. This implies a total mass change in the microscopic scale;
\item add or remove from the macroscopic scale the density lost or gained from the microscopic scale, so as to have the conservation property of the total mass of the system.
\end{itemize}
The strategy is presented in detail in Algorithm \ref{micromacro_pseudocode}.

\begin{remark}
It is important to note that, to modify the mass on a microscopic scale, it is not necessary to alter the number of particles; rather, their weight can be varied. This approach ensures the conservation of the total number of particles, thus obviating the issues associated with particle sampling.
\end{remark}

\begin{remark}

 In order for $\zeta$ to provide a correct comparison of the results obtained in the two scales and thus to allow the mass to be shifted correctly, we impose a minimum threshold value of the density that must be present in each scale. This means imposing threshold values on the function $\zeta$: 
 \begin{equation}
 \label{threshold_zeta}
     \zeta_{min}\leq \zeta(t,x) \leq \zeta_{max}, \quad \forall t\in[0,T], \forall x\in \mathbb{R}^d \text{ and with } 0<\zeta_{min}<\zeta_{max}<1.
\end{equation} 
\end{remark}

\begin{remark}
\label{movement_mass}
A potential issue that can arise from this type of micro-macro decomposition is related to the continuous and uncontrolled movement of mass between the two scales. This phenomenon is not intrinsically associated with the properties of the analyzed systems; rather, it arises from the adaptation of the microscopic and macroscopic algorithms in the presence of a mass shift. 

Also, it may be advantageous to refrain from shifting mass for small times, thereby ensuring that the algorithm initially compute results with fixed mass. In fact, it has been empirically observed that strong mass fluctuations occur for small times.
\end{remark}

\subsection{Numerical results}
\subsubsection{Constrained microscopic system}
We started by testing the algorithm on the microscopic scale using \eqref{sist_particles} to compute the evolution of the position and velocity of each particle. 
We consider the minimization problem for the Ackley function 
\begin{equation}
\label{ackley}
    \min_{x\in \mathbb{R}^d} \mathcal{F}(x):=-20 \exp\left(-0.2 \sqrt{\frac{1}{d}\sum_{i=1}^{d}(x_i)^2}\right)-\exp\left(\frac{1}{d}\sum_{i=1}^{d}cos(2 \pi (x_i))\right)+20+e
\end{equation}
that is continuous, non-differentiable and non-convex and has
its global minimizer in $x^* =(0,...,0)$. We use $N = 480$ particles for the PSO discretization in a two-dimensional search space ($d=2 $). 

\begin{figure}[H]
\centering
     \begin{subfigure}[t]{0.49\textwidth}
         \centering
         \includegraphics[width=\textwidth]{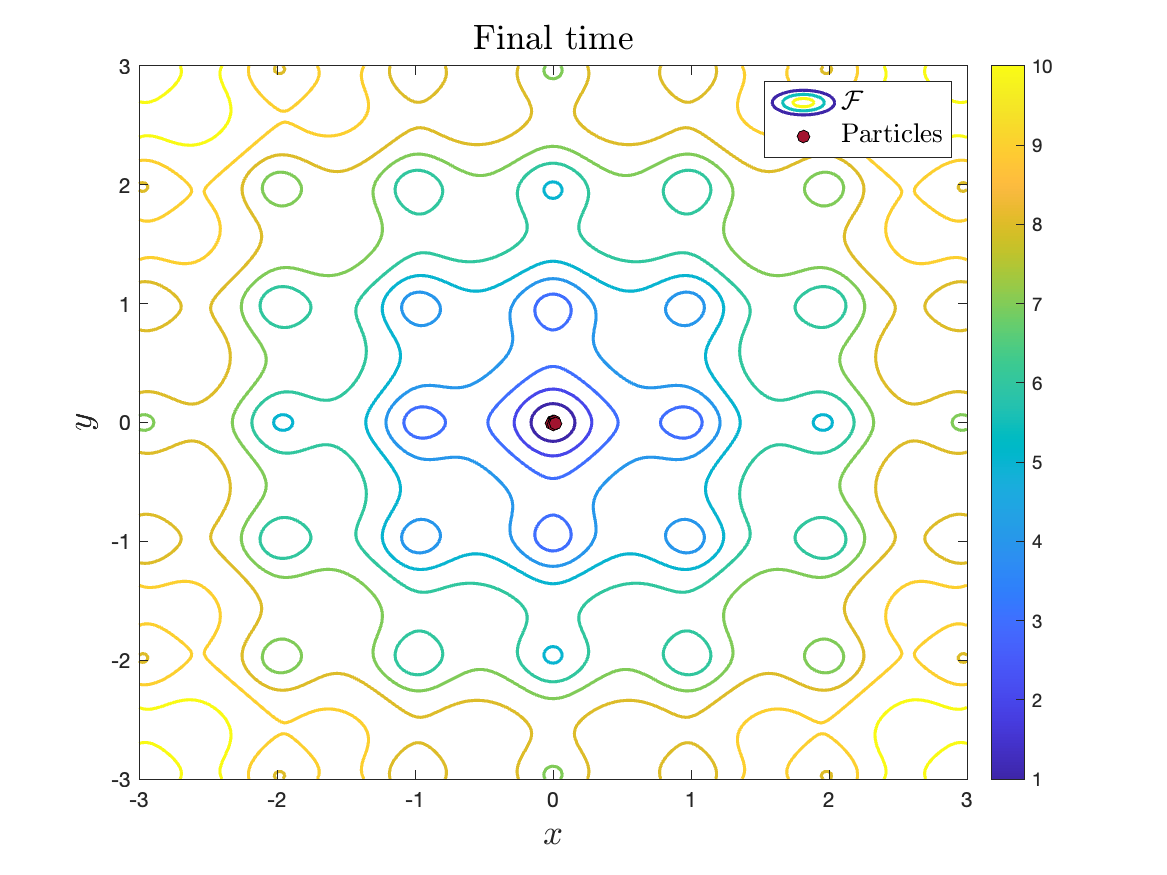}
     \caption{The values of the objective function $\mathcal{F}(x)$ \eqref{ackley} and the final positions of the particles are shown in the unconstrained problem.}
      \label{fig:nocostr}
     \end{subfigure}
     \hfill
     \begin{subfigure}[t]{0.49\textwidth}
         \centering
         \includegraphics[width=\textwidth]{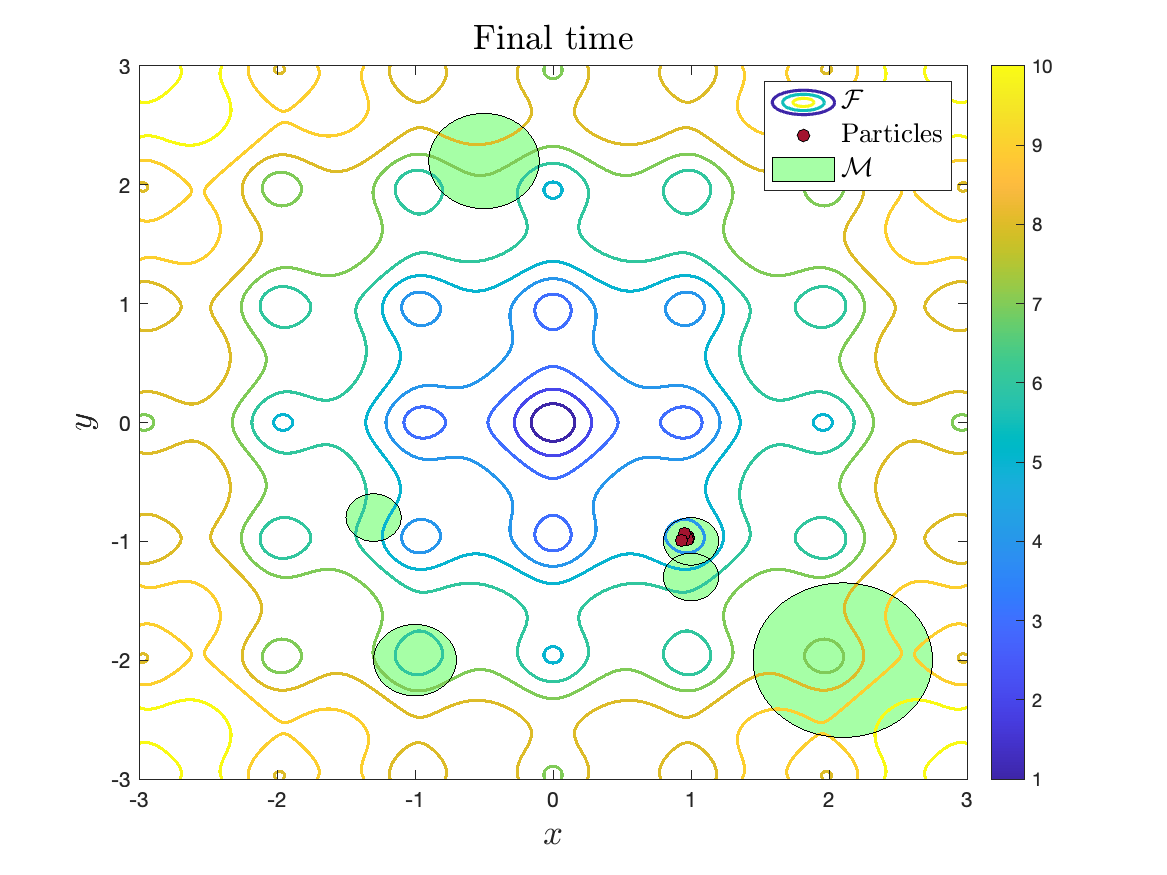}
     \caption{The values of the objective function $\mathcal{F}(x)$ \eqref{ackley}, the final positions of the particles, and the feasible set $\mathcal{M}$, which is defined as the union of the circles in green, are shown in the constrained problem.}
     \label{fig:withcostr}
     \end{subfigure}
     \hfill
     \caption {Evolution of particles in the unconstrained \ref{fig:nocostr} and constrained \ref{fig:withcostr} case, starting from the same initial conditions. The particles in the unconstrained case \ref{fig:nocostr} converge to the global minimizer $x^*$, while in the constrained case \ref{fig:withcostr} the particles converge to the global minimizer $x^*_\mathcal{M}$, that is feasible for the set $\mathcal{M}$ chosen in the simulation.}
\end{figure}

We tested the classical microscopic algorithm and the constrained algorithm on the same initial data. The following results are all obtained with the following parameters: $\gamma=0.5,$ $\lambda=1$, $\sigma=\frac{1}{\sqrt{3}}$, $\alpha=30$. The parameter values used for the penalized algorithm are as follows: $\kappa_0=5,$ $\eta_\kappa=1.1,$ $\beta_0=1,$ $ \eta_\beta=1.1$. The anisotropic diffusion is used, see \eqref{anisotropic}. In the constrained case we consider the admissible set
\begin{align*}
\label{constraint_simulations}
    \mathcal{M}=\{x\in\mathbb{R}^d|&(x_1+0.5)^2+(x_2-2.2)^2\leq 0.4 \ \cup \\&(x_1-1.3)^2+(x_2+0.8)^2\leq 0.2 \ \cup \\ &(x_1-1)^2+(x_2+1.3)^2\leq 0.1 \ \cup\\ &(x_1-1)^2+(x_2+1)^2\leq 0.1 \ \cup\\ &(x_1-2.1)^2+(x_2+2)^2\leq 0.65 \ \cup\\ &(x_1+1)^2+(x_2+2)^2\leq 0.3\}.
\end{align*}   

\begin{figure}[H]
\centering
     \begin{subfigure}[b]{0.48\textwidth}
         \centering
         \includegraphics[width=\textwidth]{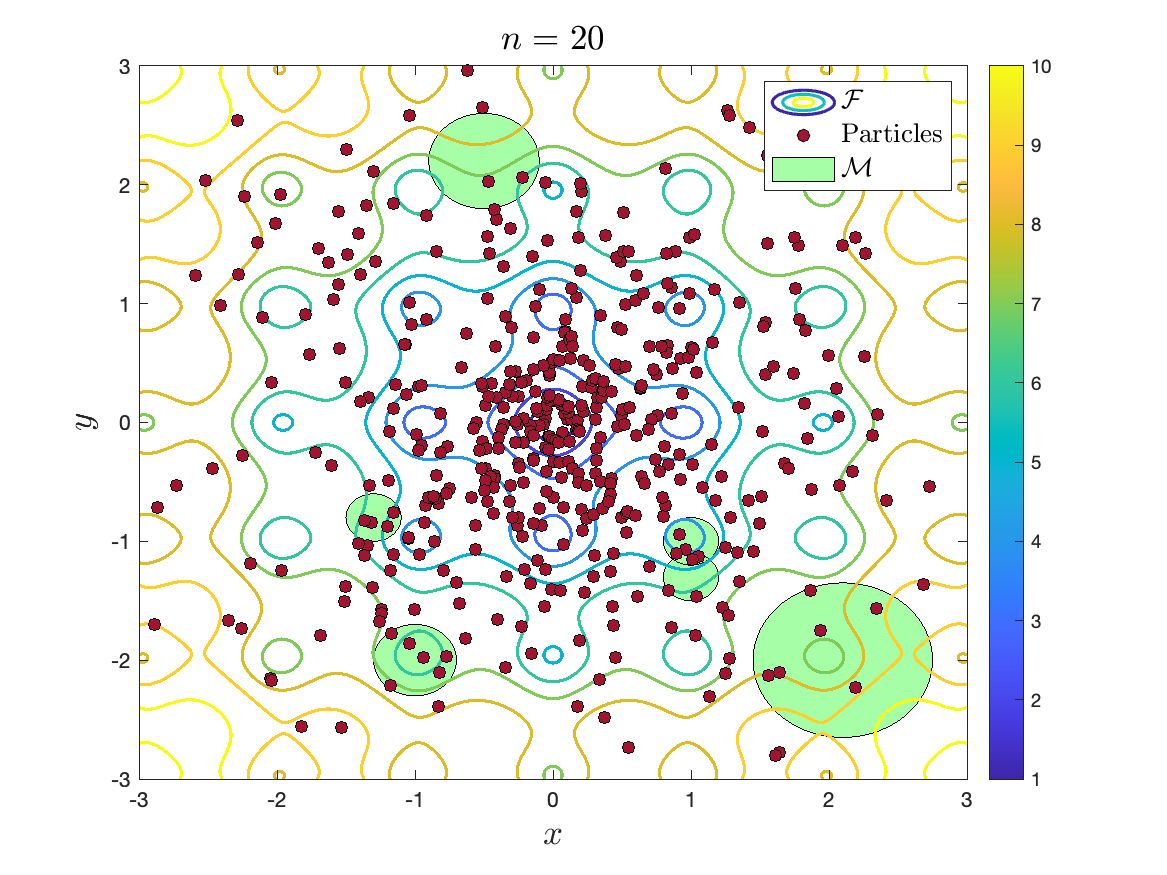}
     \caption{Positions of particles at time step $20$.}
      \label{fig:pos20}
     \end{subfigure}
     \hfill
     \begin{subfigure}[b]{0.48\textwidth}
         \centering
         \includegraphics[width=\textwidth]{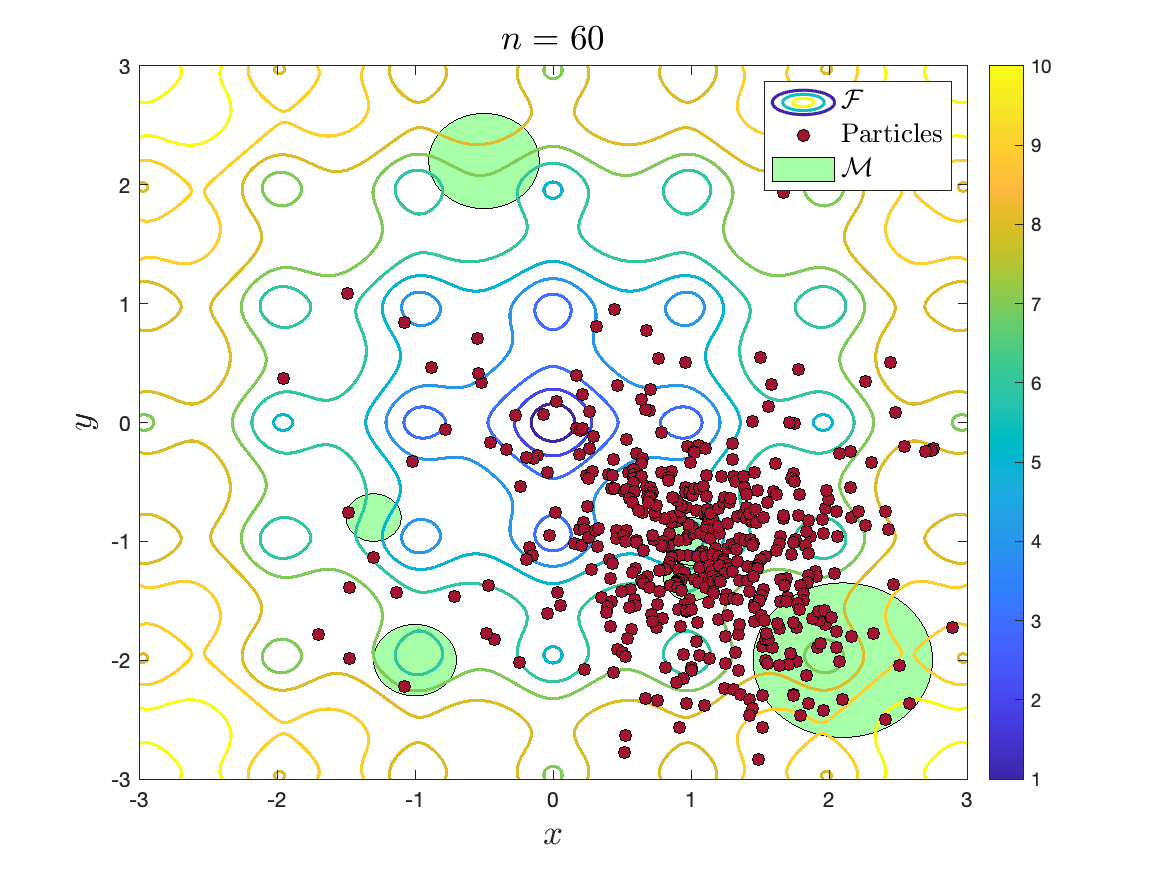}
     \caption{Positions of particles at time step $60$.}
     \label{fig:pos60}
     \end{subfigure}
     \hfill
      \begin{subfigure}[b]{0.48\textwidth}
         \centering
         \includegraphics[width=\textwidth]{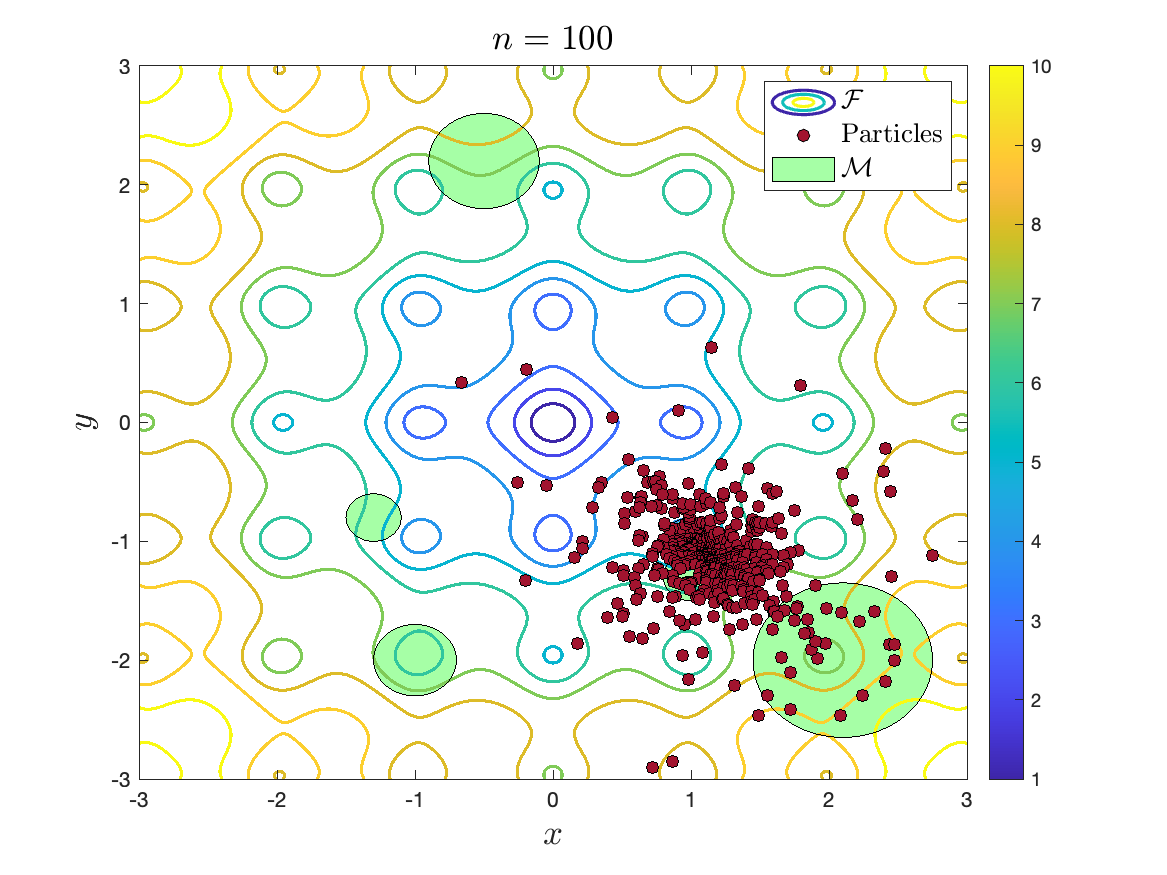}
     \caption{Positions of particles at time step $100$.}
     \label{fig:pos100}
     \end{subfigure}
     \hfill
           \begin{subfigure}[b]{0.48\textwidth}
         \centering
         \includegraphics[width=\textwidth]{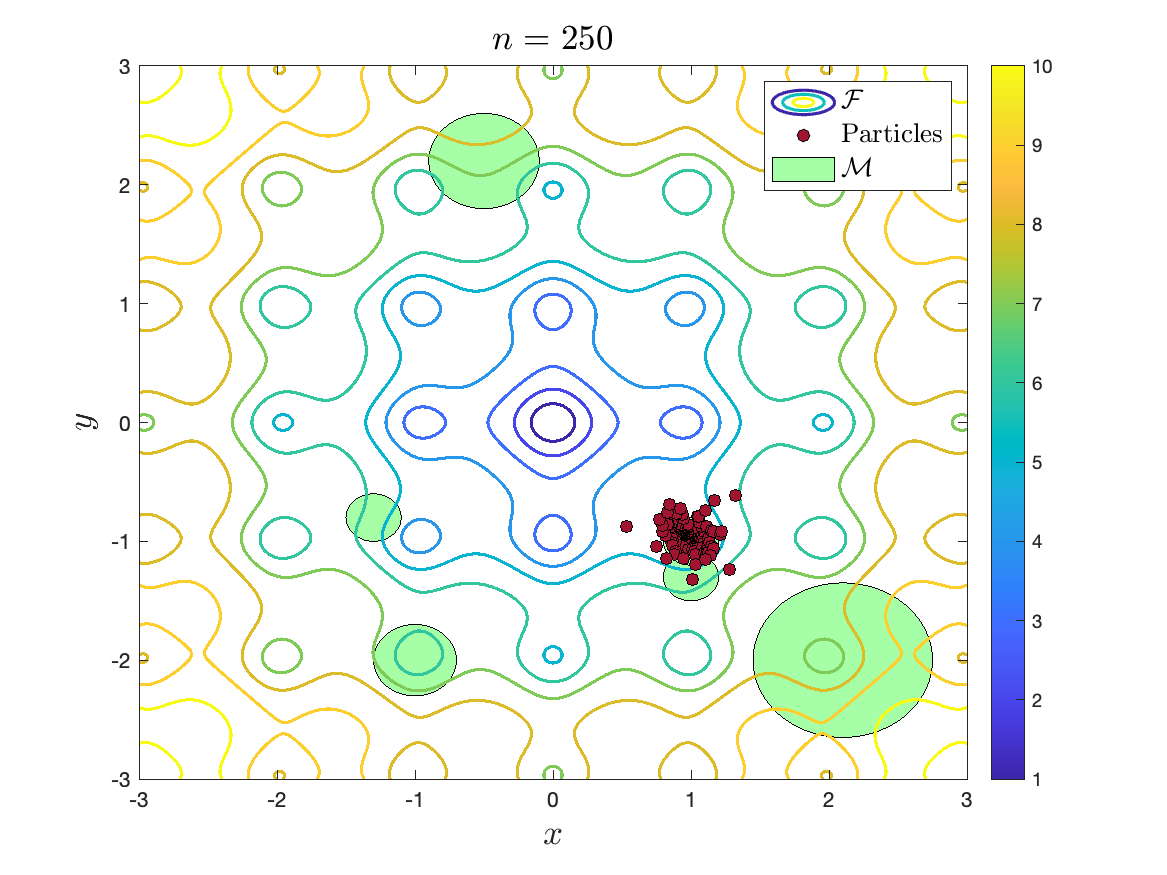}
     \caption{Positions of particles at time step $250$.}
     \label{fig:pos250}
     \end{subfigure}
     \caption {Plots show the evolution of particle positions. For small time instants, the penalty is not large enough to direct the particles toward the feasible global minimizer $x^*_\mathcal{M}$. As time evolves, the particles concentrate in $x^*_\mathcal{M}$.}
     \label{fig:posvinc}
\end{figure}

\begin{figure}[!ht]
\centering
\includegraphics[width=0.7\textwidth]{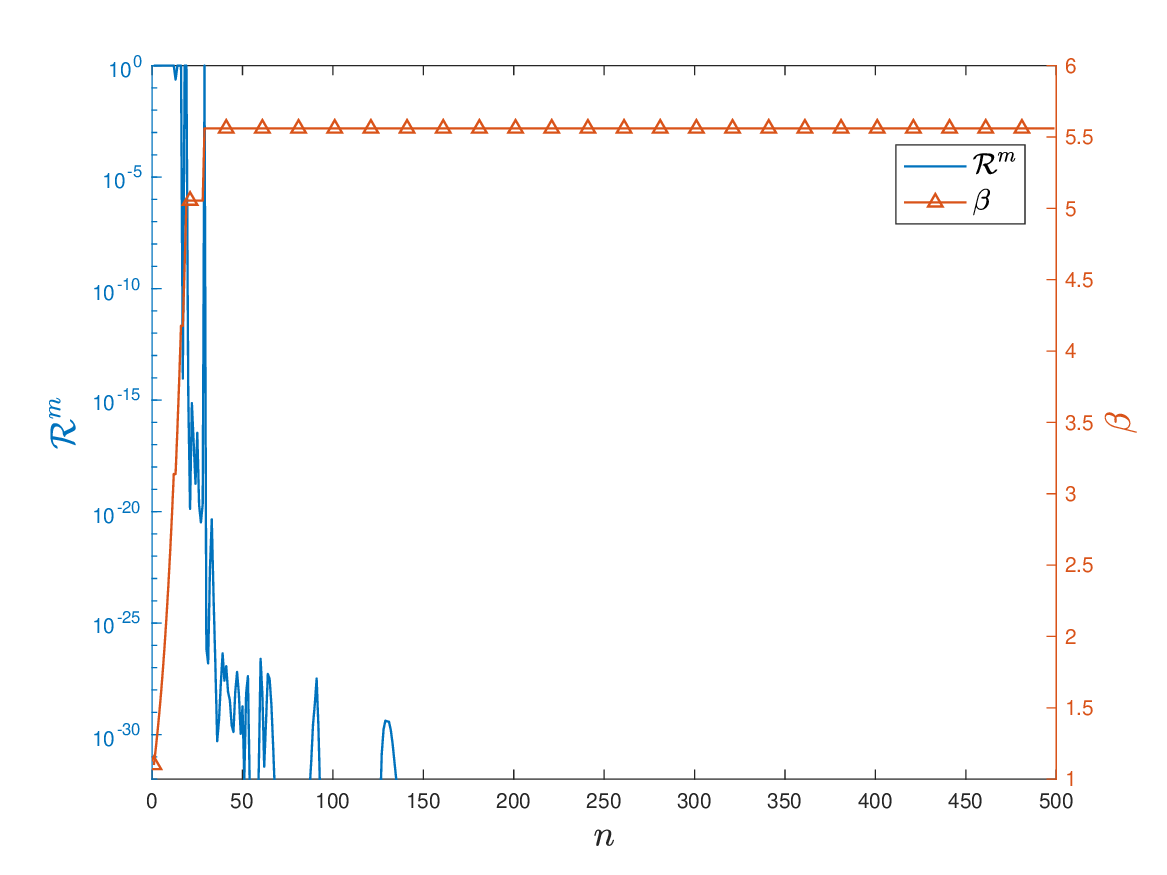}
         \caption{The evolution of two quantities is represented: the constraint violation $\mathcal{R}^m$ \eqref{R_micro} and the penalty parameter $\beta$. To effectively represent the quantity $\mathcal{R}^m$, its $y$-axis is set to a semilogarithmic scale. As we expect, the constraint violation decreases over time, while the value of $\beta$ increases.} 
         \label{fig: ev_beta}
\end{figure}

As expected, the particles at the final time $t = 400$ (shown in \Cref{fig:nocostr}), in the absence of a penalty, accumulate in the global minimizer $x^*$. In contrast, in the penalized case, all particles are attracted to the minimizer within the feasible set, $x^*_\mathcal{M}$, as illustrated in \Cref{fig:withcostr}.\\
We have plotted the positions at different intermediate time instants before convergence to the minimizer inside $\mathcal{M}$ (Figure \ref{fig:posvinc}). It is evident that initially the particles tend to converge to the global minimizer which is outside the constraint $\mathcal{M}$. This phenomenon can be attributed to the fact that the initial value of $\beta$ is substantially distant from the optimal value. Specifically, the initial value of $\beta$ is set at $\beta_0=1$ and the final value of $\beta$ attained by the algorithm is $\bar{\beta}\approx 6.1$ (see \Cref{fig: ev_beta} for a visual representation). This indicates that initially the value of the penalty term is insufficient to direct the particles within the constraint set. As time evolves, the particles reach an intermediate configuration in which they remain in proximity to the global constraint while expanding in space. Finally, as the value of $\beta$ approaches the optimal value, the particles reach the feasible global minimizer $x^*_\mathcal{M}$.\\
\Cref{fig: ev_beta} shows the evolution of the value of the parameter $\beta$ (in orange) and the constraint violation (in blue). The latter exhibits a decreasing trend, as particles tend to concentrate around the feasible global minimizer as time progresses. Consequently, the contributions of the penalty function in \eqref{R_micro} tend to be negligible. It is evident that there is a significant decrease in the value of the feasibility condition corresponding to the instant at which the value of $\beta$ increases. This marks the point in time at which particles begin to be attracted to the feasible global minimizer $x^*_\mathcal{M}$ rather than the global minimizer $x^*$ as previously.

 \begin{figure}[!ht]
\centering
\includegraphics[width=0.8\textwidth]{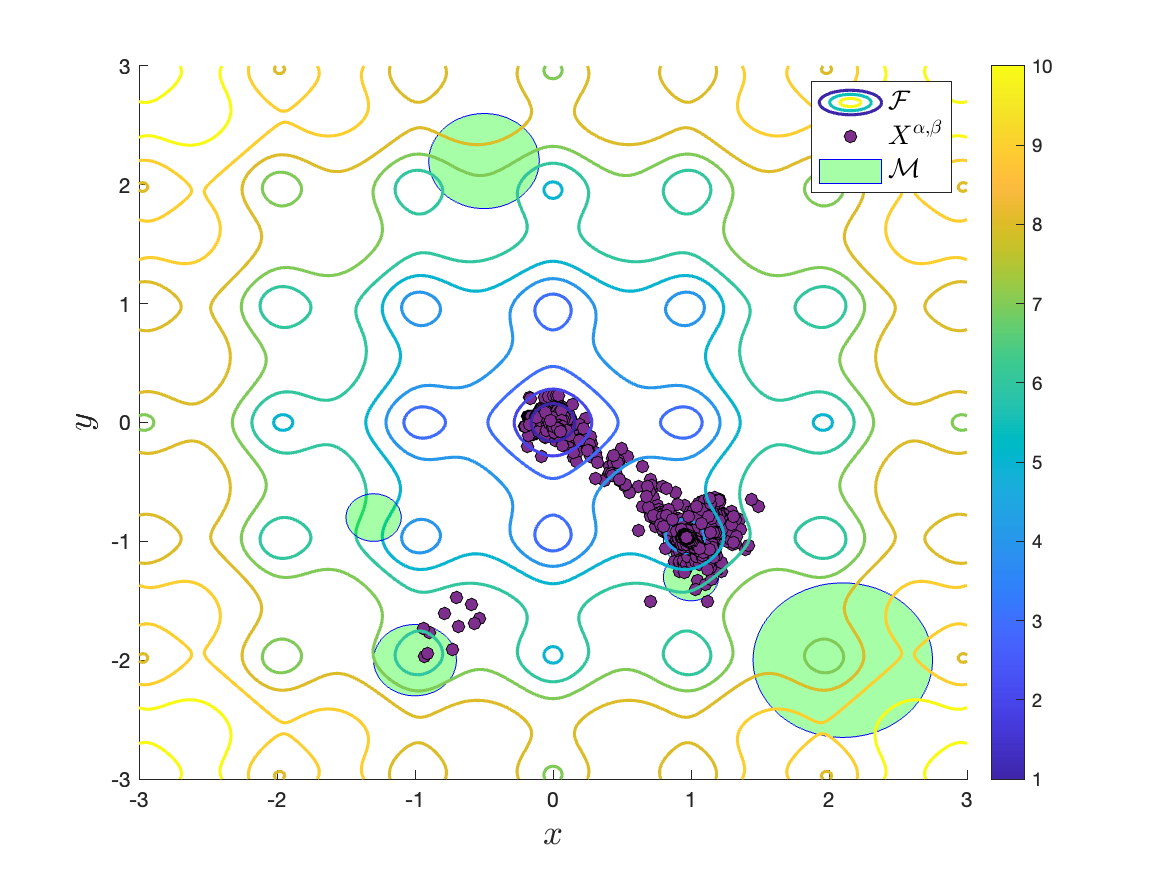}
         \caption{The plot illustrates the evolution of the consensus points obtained from 100 runs of the algorithm. It is evident in the figure, the $X^{\alpha,\beta}$ trend follows the direction of the gradient descent, whereby the algorithm initiates from the global minimizer $x^*$ of the objective function and arrives at the feasible global minimizer $x^*_\mathcal{M}$. } 
         \label{fig: gradient_behavior}
\end{figure}

We now focus on the evolution of the consensus point. In Figure \ref{fig:posvinc}, we deduced that the movement of particles from the central region of the figure toward the constraint is indicative of a better approximation of the consensus point as time evolves. To further investigate this phenomenon, we executed 100 runs of the PSO algorithm. As illustrated in \Cref{fig: gradient_behavior}, the evolution of the consensus point exhibits characteristics reminiscent of a gradient-based method, despite the implementation of the PSO method, which is inherently gradient-free. This phenomenon has previously been examined for the CBO method (which is connected to PSO via a zero inertia limit), in the context of solving unconstrained minimization problems, as outlined in \cite{riedl2023gradientneed}. \\
The presented method can be further evaluated based on the results already obtained for the PSO and CBO optimization methods. In \cite{BORGHI2023113859}, a comparison of the convergence rates  of the two methods without penalty is provided. The PSO method has been demonstrated to be less stable and reliable when dealing with difficult objective function, in contrast to the CBO method, which has been shown to be more robust.  We expect these results to remain unchanged when comparing CBO and PSO with penalization.

\subsubsection{Constrained macroscopic system}
We proceed to present the results of the simulations obtained in the macroscopic case. We once again consider the minimization problem of the Ackley function (see \eqref{ackley}) on a one-dimensional space with a grid of 401 points in space and with a constraint given by
\begin{equation*}
    \mathcal{M}=\{x\in\mathbb{R}|-1.8\leq x\leq -1.6 \cup -1.2\leq x\leq -0.8 \cup 1.1\leq x\leq 1.3 \cup 1.7\leq x\leq 1.9 \}.
\end{equation*}
The same set of parameters employed in the microscopic case is also utilized here. The CFL condition chosen for all macroscopic simulations is equal to 0.8. In the simulations, we use the second-order macroscopic system by applying the Lax-Friedrichs numerical scheme \eqref{eqn: LxF}.

\begin{figure}[H]
\centering
     \begin{subfigure}[b]{0.55\textwidth}
         \centering
         \includegraphics[width=\textwidth]{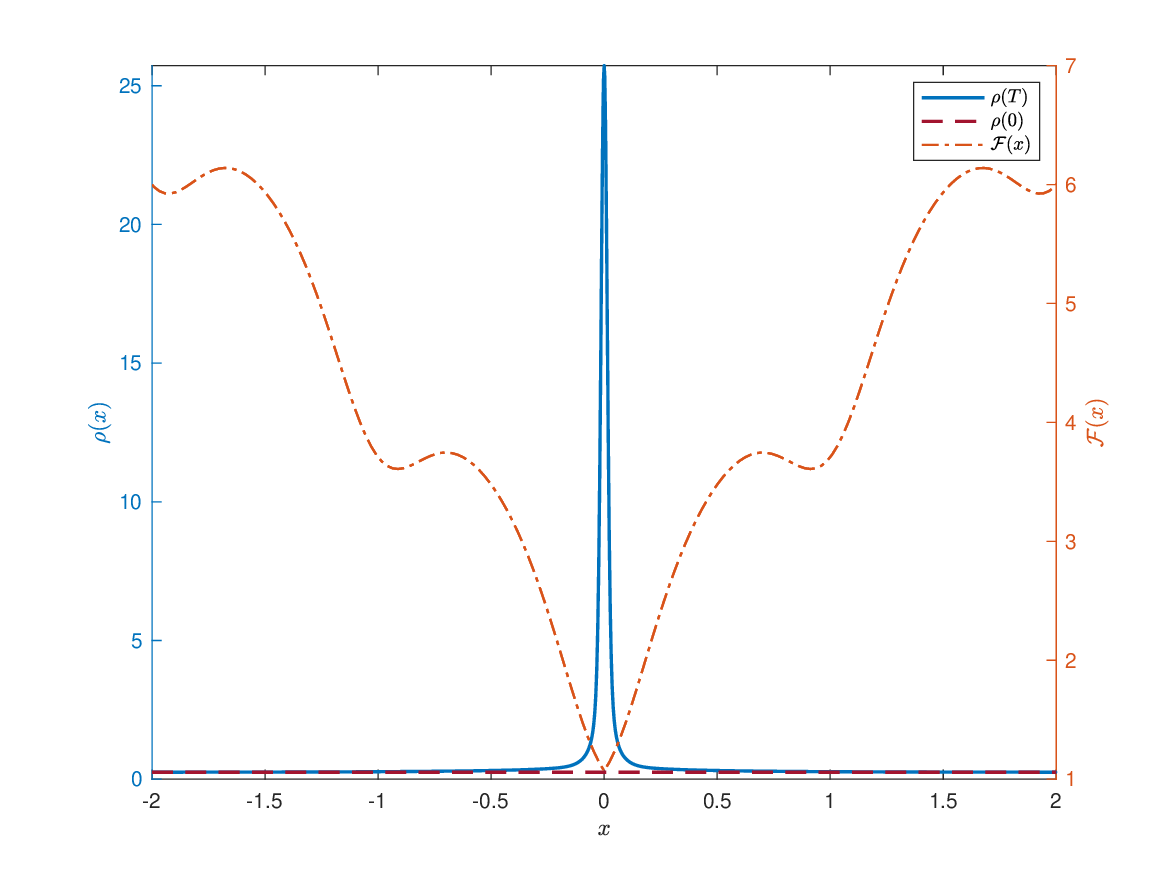}
\caption{In the unconstrained minimization problem, we see that the density $\rho$ converges to a Dirac delta centered in $0$, which is the global minimizer $x^*$.}
\label{fig:macro1}
     \end{subfigure}
     \hfill
     \begin{subfigure}[t]{0.49\textwidth}
         \centering
         \includegraphics[width=\textwidth]{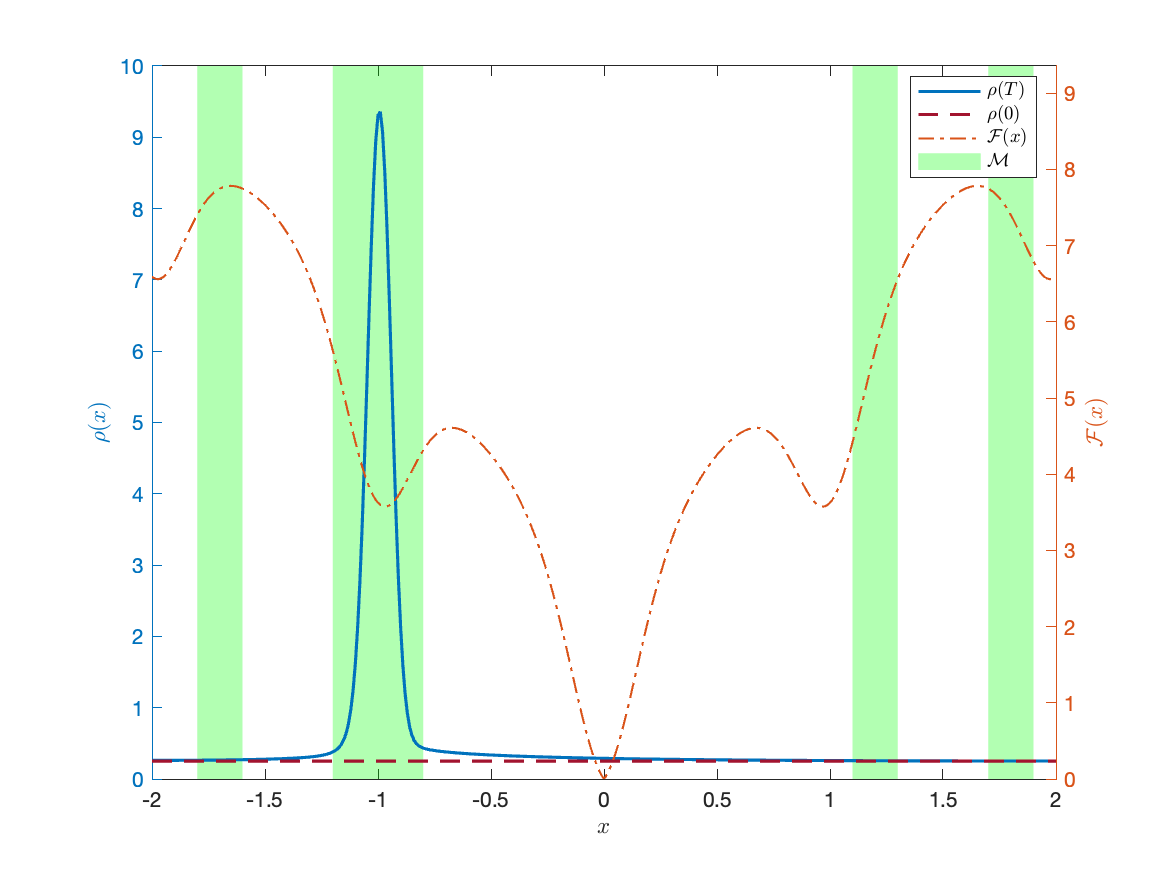}
     \caption{In the constrained minimization problem, we see that the density $\rho$ converges to a Dirac delta centered in $-1$, which is the feasible global minimizer $x^*_\mathcal{M}$.}
     \label{fig:macro2}
     \end{subfigure}
     \hfill
     \begin{subfigure}[t]{0.49\textwidth}
         \centering
        \includegraphics[width=\textwidth]{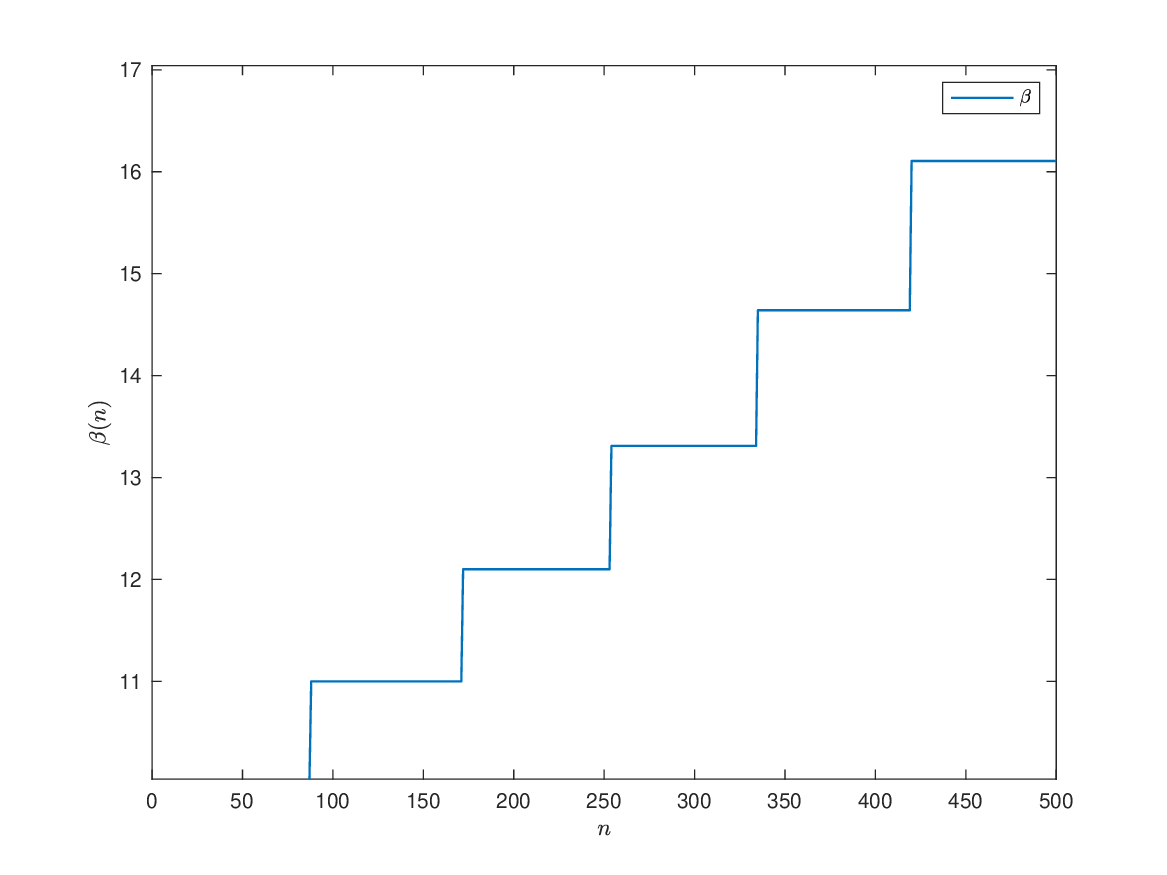}
     \caption{Evolution of the value of the penalty strength $\beta$ that has an increasing trend, as expected.}
     \label{fig:beta_evol_macro}
     \end{subfigure}
     \caption{Plots obtained by applying the PSO method to constrained and unconstrained minimization problems in which the objective function is the Ackley function.}
\end{figure}

\Cref{fig:macro1} and \Cref{fig:macro2} show the initial and final densities obtained by the algorithm in the unconstrained and constrained cases. As previously highlighted, in the unconstrained scenario depicted in Figure \ref{fig:macro1}, the mass is concentrated at the global minimizer $x^*$, in fact the final density tends to a Dirac delta centered at $x^*$. Instead in the constrained case \ref{fig:macro2}, the final density tends to a Dirac delta centered at the minimizer within the feasible set $\mathcal{M}$. \\
The time evolution of the penalty parameter $\beta$ is shown in \Cref{fig:beta_evol_macro}. The value of the parameter increases over time, which is a feature of the algorithm. This phenomenon is analogous to that observed in the microscopic algorithm (\Cref{fig: ev_beta}). It is important to note that the penalized problem we are solving yields the same result as the constrained problem we are trying to solve.

\subsubsection{Micro-macro decomposition}
\begin{figure}[H]
\centering
     \begin{subfigure}[b]{0.6\textwidth}
         \centering
         \includegraphics[width=\textwidth]{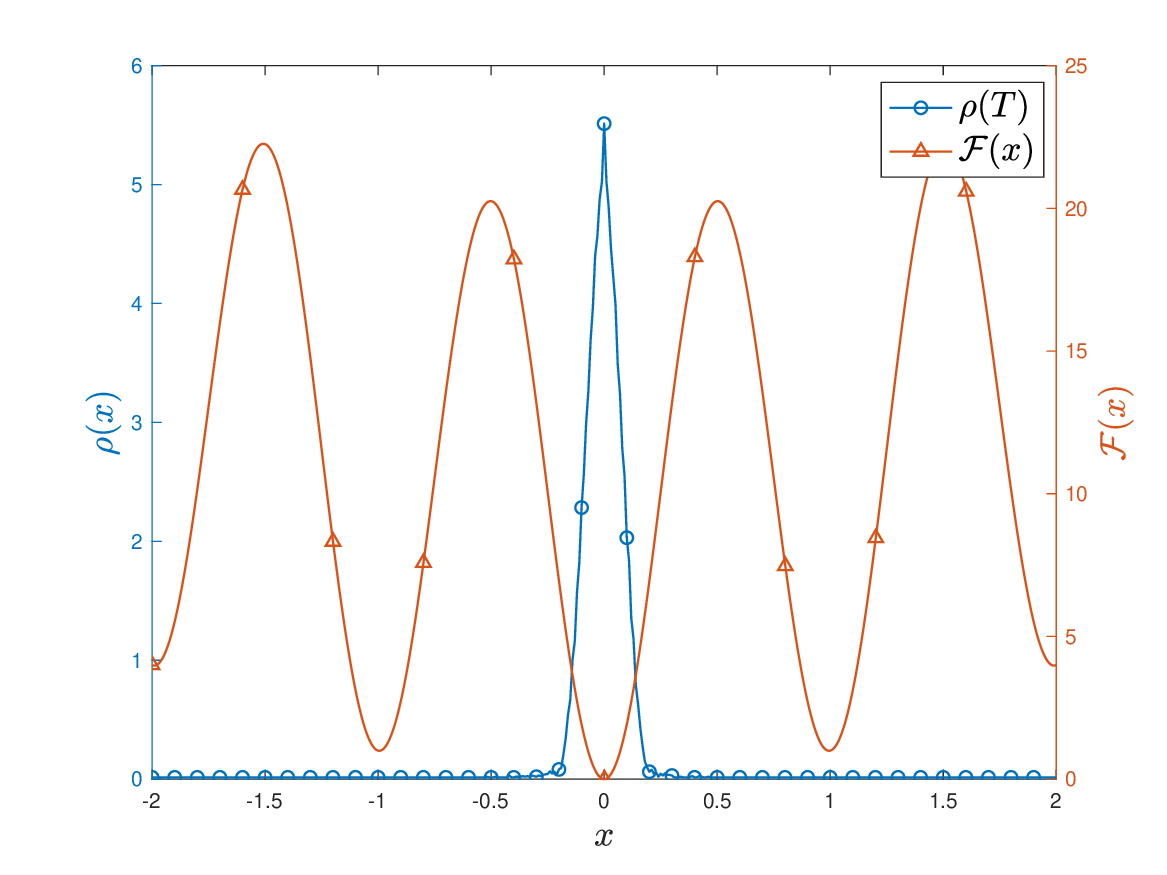}
         \caption{The density $\rho$ converges to a Dirac delta centered in $0$, which is the global minimizer.}
         \label{fig:micromacro1}
     \end{subfigure}
     \hfill
     \begin{subfigure}[t]{0.48\textwidth}
         \centering
         \includegraphics[width=\textwidth]{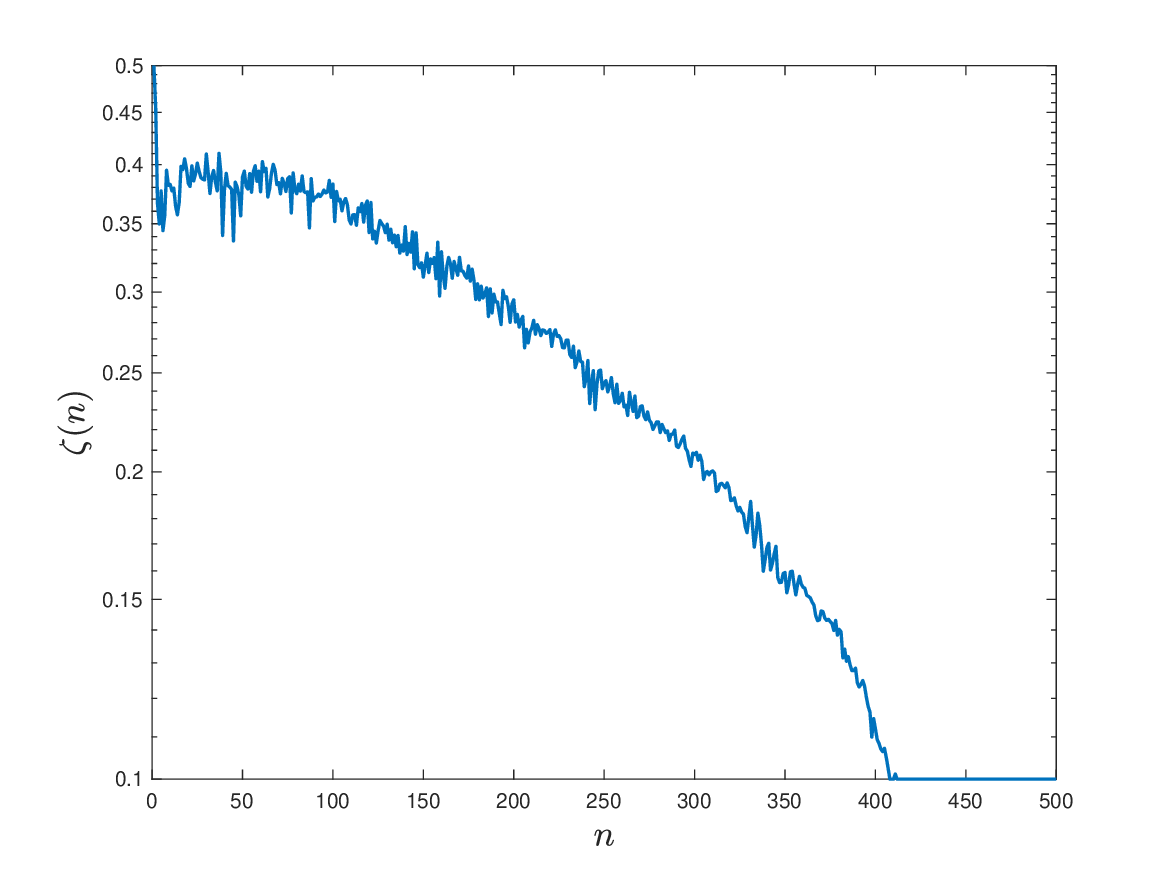}
         \caption{The value of $\zeta$ decreases, meaning that the microscopic contribution is dampened as time evolves. }\label{fig:micromacro2}
     \end{subfigure}
     \hfill
     \begin{subfigure}[t]{0.48\textwidth}
         \centering
         \includegraphics[width=\textwidth]{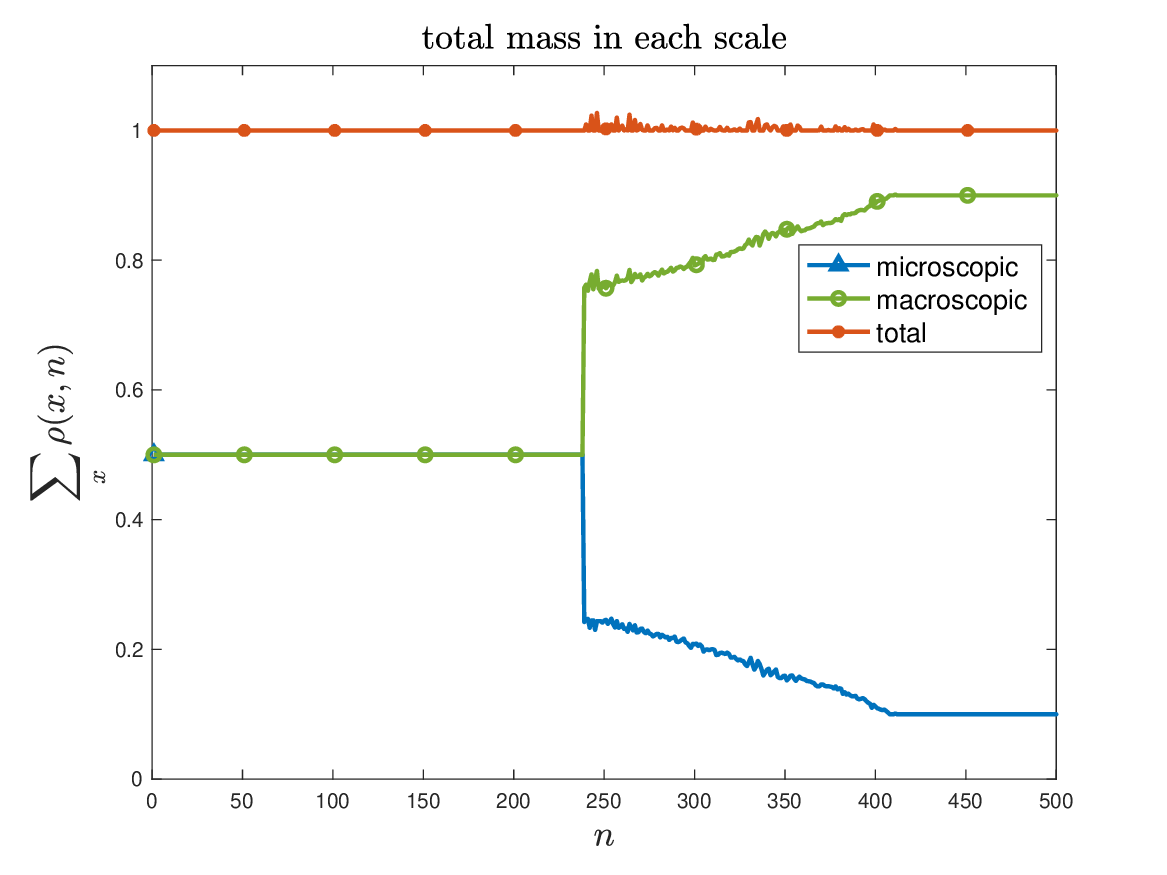}
         \caption{We see that up to time step $240$ there is no movement of mass between the two scales. At the limit, mass tends to move toward the macroscopic scale.}
         \label{fig:micromacro33}
     \label{fig:mass_evol_micromacro}
     \end{subfigure}
     \caption{Plots obtained by applying micro-macro decomposition to solve an unconstrained minimization problem with objective function the Rastrigin function \eqref{eqn: Rastrigin_function}.}
\end{figure}

We consider the minimization problem for the Rastrigin function 
\begin{equation}
\label{eqn: Rastrigin_function}
     \min_{x\in \mathbb{R}^d} \mathcal{F}(x):=10d+\sum_{i=1}^{N} [x_i^2-10 cos (2\pi x_i)]
\end{equation}
that is continuous, differentiable and non-convex, has many local minima and a global minimizer in $x^* =(0,...,0)$. We used $N = 480$ particles for PSO discretization for a one-dimensional search space ($d=1$).
In all the simulations following that utilize micro-macro decomposition, the threshold values \eqref{threshold_zeta} are set to $\zeta_{min}=0.1$ and $\zeta_{max}=0.9$.
The value of $\zeta$ at the initial instant is fixed at 0.5, ensuring equivalent initial density in microscopic and macroscopic systems. As mentioned in Remark \ref{movement_mass}, to improve the stability of the algorithm, we avoid mass shifts from the first time steps.  In the simulations, we set $t^*=240$ as the first instant at which we modify the value of $\zeta$ and subsequently move the mass between the microscopic and macroscopic scales. The choice of $t^*$ has a direct impact on the speed of convergence.\\
\Cref{fig:micromacro1} shows the final total density of the system studied. Close to the global minimizer $x^*$, one can distinctly see the contribution that the microscopic scale makes to the profile of the density distribution.
\Cref{fig:micromacro2} shows the value of $\zeta$, which, as defined, represents when the results obtained in the two scales differ from each other. It can be seen that this value tends to zero, as in the limit the results obtained from the two systems tend to the same final configuration.
In fact, by definition \eqref{def: zeta}, $\zeta$ measures the deviation of the particles from the chosen closure function for the macroscopic system, which in our work is the Maxwellian distribution \eqref{eqn: macrosyst_closure}. The fact that $\zeta$ tends to zero means that the particles tend to the chosen closure function.\\
In \Cref{fig:micromacro33} there is a plot of the mass present in each scale at each time instant. As previously mentioned, the initial mass on the two scales is equivalent up to the constant value, designated as $t^*=240$. It should be noted that this choice was determined empirically to ensure a nearly monotonic movement of the mass between the two scales, thereby accelerating the achievement of convergence. Actually, a shift in mass towards the macroscopic scale is observed, accompanied by a decrease in the microscopic scale. This outcome shows that, as time progresses, particle interactions become less significant in the problem under analysis. This study may be extended to address a broader spectrum of minimization problems. In high-dimensional minimization problems, we expect that mass will not always concentrate at the macroscopic scale. Due to the inherent complexity and challenges posed by high-dimensional macroscopic behavior, mass may tend to concentrate at the microscopic scale.

\subsubsection{Constrained micro-macro decomposition}
It has been shown that micro-macro decomposition can be used to address minimization problems in which the domain of the objective function is the entire space, denoted by $\mathbb{R}^d$. Furthermore, its application in constrained penalty problems, which are addressed by the adaptive penalized strategy outlined in Section \ref{sec:micro_macro_decomposition}, is also possible. At the algorithmic level, at each time step, the strategy consists of solving the penalized problem on the microscopic and macroscopic scales. Then we proceed with the comparison of the results in the two scales following Algorithm \ref{micromacro_pseudocode}. Subsequently, the procedure is repeated for the next time step.\\
We present a simulation in which the chosen objective function is the Rastrigin function \eqref{eqn: Rastrigin_function} in one-dimensional space and the chosen constraint is
\begin{equation*}
    \mathcal{M}=\{x\in\mathbb{R}|x<-0.5\}.
\end{equation*}
All parameters are chosen as in previous simulations.\\
\begin{figure}[H]
\centering
     \begin{subfigure}[b]{0.6\textwidth}
         \centering
         \includegraphics[width=\textwidth]{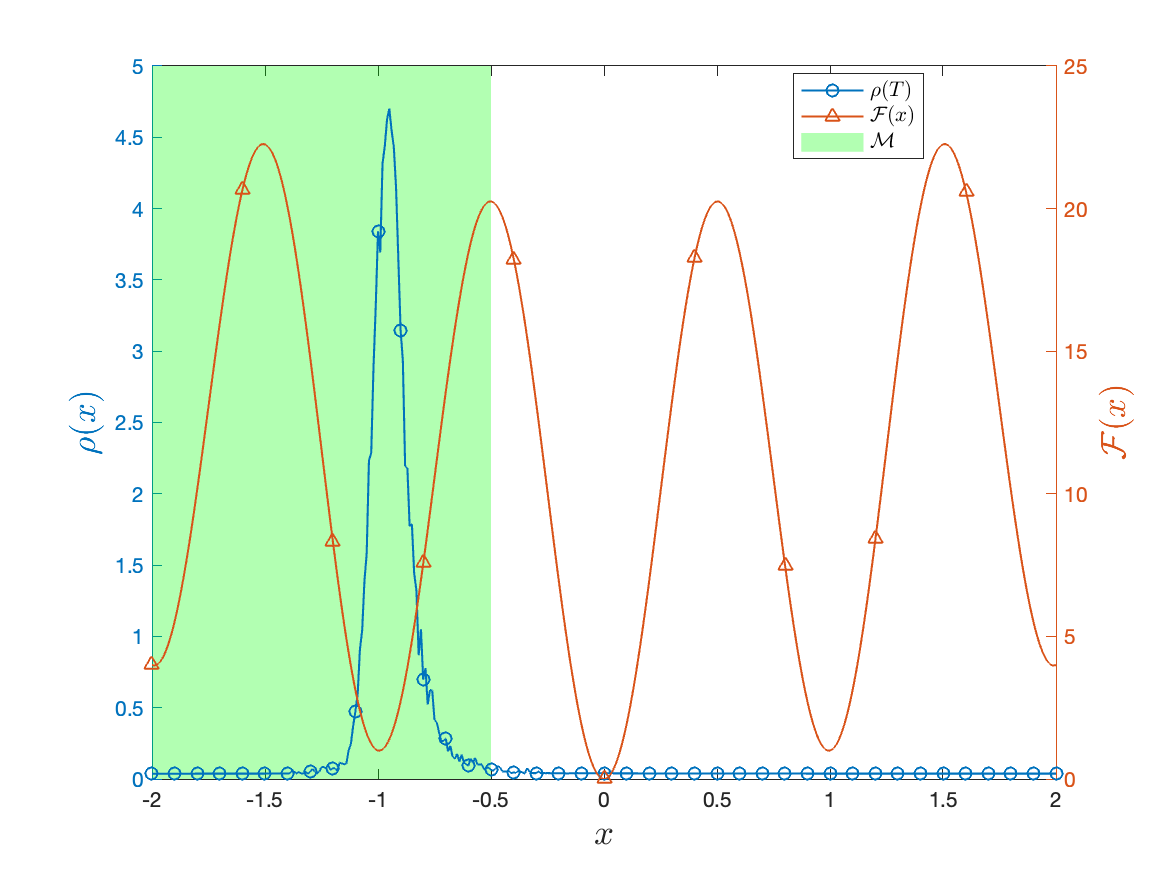}
         \caption{The density $\rho$ converges to a Dirac delta centered in $-1$, which is the feasible global minimizer.}
         \label{fig:micromacro3}
     \end{subfigure}
     \hfill
     \begin{subfigure}[t]{0.48\textwidth}
         \centering
         \includegraphics[width=\textwidth]{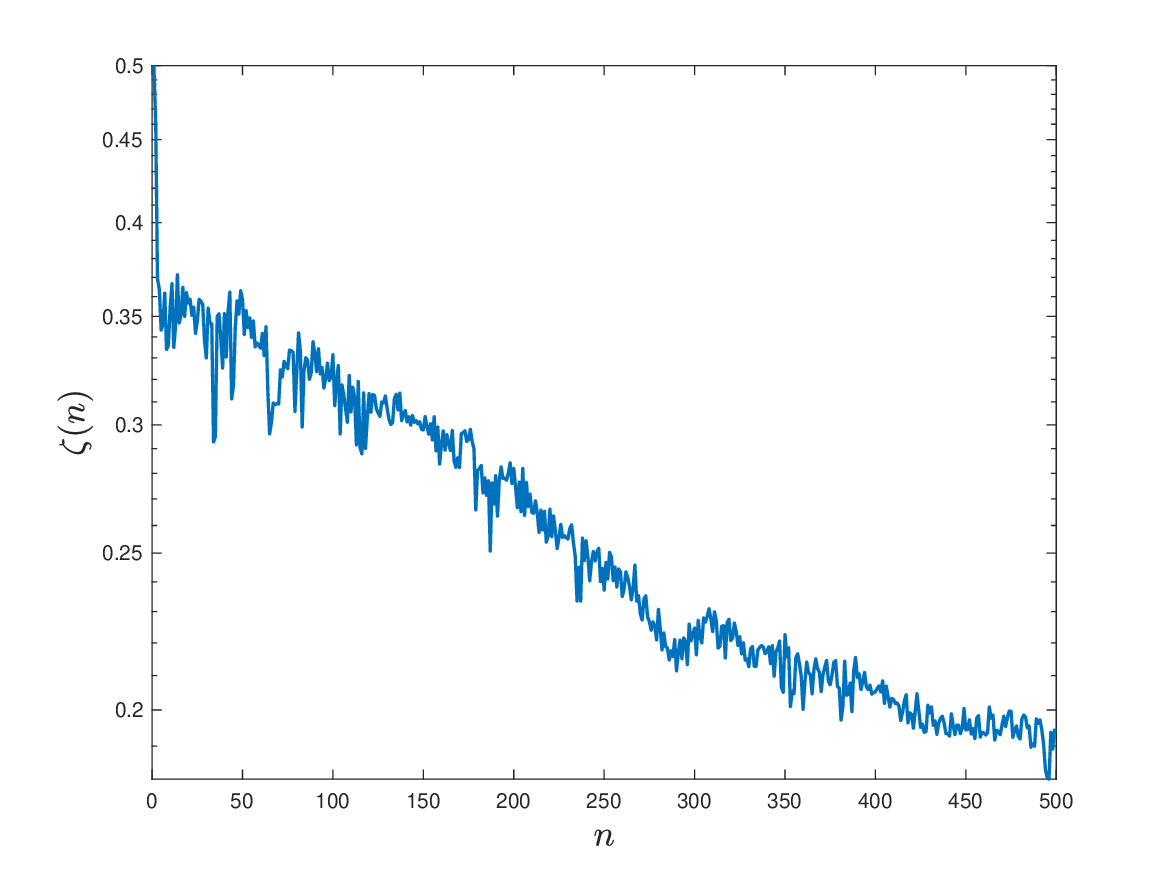}
         \caption{The value of $\zeta$ decreases, meaning that the microscopic contribution is dampened as time evolves because it is multiplied by $\zeta$. }\label{fig:micromacro4}
     \end{subfigure}
     \hfill
     \begin{subfigure}[t]{0.48\textwidth}
         \centering
         \includegraphics[width=\textwidth]{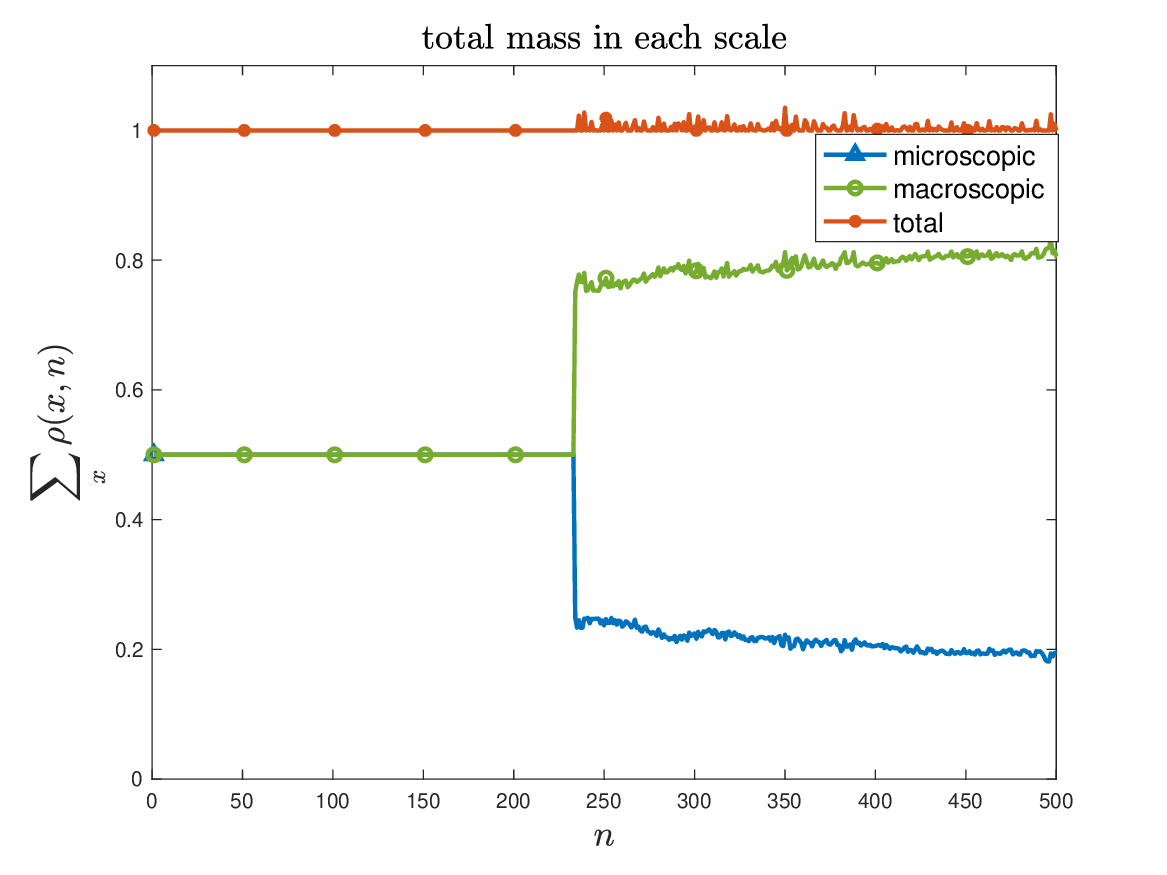}
         \caption{As explained in Remark \ref{movement_mass}, we see that up to time step $240$ there is no mass movement between the two scales. At the limit, mass tends to move toward the macroscopic scale.}
         \label{fig:micromacro5}
     \label{fig:mass_evol_micromacro_constr}
     \end{subfigure}
     \caption{Plots obtained by applying micro-macro decomposition to solve a constrained minimization problem with objective function the Rastrigin function \eqref{eqn: Rastrigin_function}.}
\end{figure}

In Figure \ref{fig:micromacro3} we see that the final density, which is given by the microscopic and macroscopic contributions, converges to a Dirac delta centered at $x=-1$, which is exactly the feasible global minimizer. Again, the mass tends to shift from the microscopic to the macroscopic scale (see Figure \ref{fig:micromacro4}). This phenomenon indicates that, as time progresses, particle interactions become less significant for the purposes of convergence.

\section{Conclusion}
\label{sec:conclusion}
In this paper we conducted an analysis of the PSO method \cite{kennedy,kennedy2,shi, huang}, which is a metaheuristic method for solving minimization problems. A novel extension to the method is presented that can be applied in the case of constrained minimization problems.  This extension was introduced for the CBO model \cite{borghi2021constrained}. The strategy consists of working with an unconstrained minimization problem with a novel objective function that incorporates a penalty term. This penalty is in addition to the objective function of the initial minimization problem, and it is adaptively computed to reach its optimal value. The penalized adaptive algorithm was implemented at both the microscopic and macroscopic levels. A subsequent extension of this research would involve applying the algorithm to a macroscopic system derived from an alternative choice of closure that satisfies the conditions presented in the paper.\\
Additionally, the micro-macro decomposition is presented in this work. It is a strategy designed to understand the different relevance of particle interactions as time evolves. Inspired by the works \cite{degond, lemoumehats, dimarcopareschi}, this approach involves working in parallel on the microscopic and macroscopic levels. In each scale, a mass quantity is calculated at each time step by comparing the results obtained at the previous time step across all scales. We have shown that the micro-macro decomposition can also be applied to the constrained case.\\ Numerous results and numerical applications are presented to validate the good performance of the proposed strategies.

\section*{Acknowledgements}
The work of S.V. is supported by the European Union’s Horizon Europe research and innovation program under the Marie Sklodowska-Curie Doctoral Network Datahyking (Grant No. 101072546).
S.V. is a member of the INdAM Research National Group of Scientific Computing
(INdAM-GNCS). M.H. gratefully acknowledges support from the DFG through the projects HE5386/33-1 Control of Interacting Particle Systems, and Their Mean-Field, and Fluid-Dynamic Limits (560288187), and HE5386/34-1 Partikelmethoden für unendlich dimensionale Optimierung (561130572). \\
S.V. would like to thank René Pinnau and Giuseppe Visconti for the stimulating discussions during the preparation of this manuscript.


\end{document}